\documentclass[12pt]{amsart}
\usepackage{graphicx,amssymb,amsmath,amsfonts,amscd}
\newtheorem{thm}{Theorem}[section]
\newtheorem{cor}[thm]{Corollary}
\newtheorem{lem}[thm]{Lemma}
\newtheorem{prop}[thm]{Proposition}

\newcommand{\Trace}{\mathrm{Trace}}

\newcommand{\isomto}{\overset{\sim}{\rightarrow}}
\newcommand{\Sym}{\mathrm{Sym}}
\newcommand{\Fq}{\mathbb F_q}

\newcommand{\QQ}{\bar{\mathbb Q}_\ell}
\newcommand{\R}{\mathrm R}
\newcommand{\Zp}{\mathbb Z_p}
\newcommand{\Z}{\mathbb Z}
\newcommand{\Gal}{\mathrm {Gal}}
\newcommand{\AAA}{\mathbb A}
\newcommand{\PP}{\mathbb P}
\newcommand{\GG}{\mathbb G}

\newcommand{\HH}{\mathrm H}

\newcommand{\LL}{{\mathcal L}}

\newcommand{\Tr}{\mathrm {Tr}}

\newcommand{\Dbc}{{\mathcal D}^b_c(\AAA^1_k,\QQ)}

\newcommand{\xx}{\mathbf x}
\newcommand{\yy}{\mathbf y}
\newcommand{\vv}{\mathbf v}
\newcommand{\ww}{\mathbf w}

\title{Moment Zeta Functions For Toric Calabi-Yau Hypersurfaces}
\author{Antonio Rojas-Leon}
\address{Department of Mathematics, University of California, Irvine, CA 92697-3875, USA}
\email{arojasl@math.uci.edu}
\author{Daqing Wan}
\address{Department of Mathematics, University of California, Irvine, CA 92697-3875, USA}
\email{dwan@math.uci.edu}

\begin{document}
\maketitle
\section{Introduction}

Let $n\geq 2$ be a positive integer. We consider the
following family
$$X_{\lambda}: x_1 +\cdots +x_n + {1\over x_1\cdots x_n}
=\lambda$$ of $(n-1)$-dimensional toric Calabi-Yau hypersurfaces
in ${\mathbb G}_m^n$ parameterized by $\lambda\in {\mathbb A}^1$.
Let ${\mathbb P}_{\Delta}$ be the projective toric variety
associated to the Newton polytope of the above Laurent polynomial.
The projective closure $Y_{\lambda}$ of $X_{\lambda}$ in ${\mathbb
P}_{\Delta}$ is simply the quotient by $G=(\Z /(n+1)\Z)^{n-1}$ of
the following Dwork family of projective Calabi-Yau hypersurfaces in ${\mathbb
P}^n$:
$$W_{\lambda}: x_0^{n+1}+\cdots +x_n^{n+1} = \lambda x_0\cdots
x_n.$$ The crepant resolution of the family $Y_{\lambda}$ is the
mirror family of $W_{\lambda}$.

Let $\Fq$ be a finite field of $q$ elements with characteristic
$p$. In this paper, we are interested in the moment zeta function
which measures the arithmetic variation of the zeta function of
$X_{\lambda}$ over $\Fq$ as $\lambda$ varies in $\Fq$. The moment
zeta function grew out of the second author's study
\cite{W20}\cite{W21}\cite{W2} of Dwork's unit root conjecture. Its
general properties were studied in Fu-Wan \cite{fw0} and Wan
\cite{W0}\cite{W3}. Note that the zeta function of $Y_{\lambda}$
differs from the zeta function of $X_{\lambda}$ by some trivial
factors.

The zeta function of the Dwork family $W_{\lambda}$ over finite fields
had been studied extensively in the literature, first by Dwork \cite{Dw0} and
Katz \cite{ka0}, and more recently in connection with arithmetic mirror
symmetry by Candelas, de la Ossa and Rodriques-Villegas \cite{Can1}\cite{Can2},
and by Wan \cite{W4}\cite{W5} and Fu-Wan \cite{fw3}. By \cite{W4}\cite{W5},
the zeta function of $X_{\lambda}$ is the most primitive
piece of the zeta function of $W_{\lambda}$.
Thus, we shall restrict ourself to the family $X_{\lambda}$.
The Hasse-Weil zeta function (but not its higher moment zeta function
which would seem to be too hard at the moment) in a similar number
field example is studied in a recent paper by Harris,
Shepherd-Barron and Taylor \cite{Ta}.

More precisely, for a positive integer $d$, let $N_d(k)$ denote
the number of points on the family $X_{\lambda}$ such that $x_i\in
{\mathbb F}_{q^{dk}}$ for all $1\leq i\leq n$ and $\lambda \in
{\mathbb F}_{q^k}$. The $d$-th moment zeta function of the
morphism $X_{\lambda}\rightarrow \lambda \in {\mathbb A}^1$
 is defined to be
$$Z_d({\mathbb A}^1, X_{\lambda}) =\exp(\sum_{k=1}^{\infty} {N_d(k)\over
k}T^k)\in 1+T{\mathbb Z}[[T]].$$ This sequence $Z_d({\mathbb A}^1,
X_{\lambda})$ ($d=1,2,\cdots)$ of power series gives a simple diophantine
reformulation on the arithmetic variation of the zeta function of
the family $X_{\lambda}$. It is a rational function in $T$ for
each $d$. In the special case $n=2$, $X_{\lambda}$ is a family of
elliptic curves and the moment zeta function $Z_d({\mathbb A}^1,
X_{\lambda})$ is closely related to arithmetic of modular forms.
In general, Dwork's unit root zeta functions \cite{Dw} attached to
this family are the $p$-adic limits of this sequence of moment
zeta functions. They are thus infinite $p$-adic moment zeta
functions in some sense. Our aim of this paper is to give a
precise study of this sequence $Z_d({\mathbb A}^1, X_{\lambda})$
and their $p$-adic variation as $d$ varies $p$-adically. One main
consequence of our results is a complete determination of the
purity decomposition and the trivial factors for the moment zeta
function $Z_d({\mathbb A}^1, X_{\lambda})$ for all $d$, all $n$
and all $p$ not dividing $n+1$. This provides the first higher dimensional
example for which all higher moment zeta functions are determined.

\begin{thm} Assume that $p$ does not divide $n+1$. Then, the
$d$-th moment zeta function has the following factorization
$$Z_d({\mathbb A}^1, X_{\lambda})^{(-1)^{n-1}} =P_d(T)Q_d(T),$$
where $Q_d(T)$ is the trivial factor given explicitly by
$$\frac{
(1-q^{\frac{d(n-1)}{2}}T)^{\frac{1+(-1)^{d+n}}{2}}}
{(1-q^{\frac{d(n-1)}{2}+1}T)^{\frac{-(-1)^n-(-1)^{n+d}}{2}}}
\prod_{k=0}^{[\frac{n-2}{2}]}\frac{1-q^{dk}T}{1-q^{dk+1}T}
\prod_{i=0}^{n-1}(1-q^{di+1}T)^{(-1)^{i+1}{n \choose i+1}},$$ and
$P_d(T)$ is the non-trivial factor which has the form
$$P_d(T)=\prod_{a+b=d, 0\leq b\leq n}P_{a,b}(T)^{(-1)^{b-1}(b-1)},$$
where each $P_{a,b}(T)$ is a polynomial in $1+T{\Bbb Z}[T]$, pure
of weight $d(n-1)+1$, whose degree is given explicitly in Theorem
3.10.
\end{thm}

\begin{cor}Assume that $p$ does not divide $n+1$. Let $N_d(k)$ denote
the number of points on the family $X_{\lambda}$ such that $x_i\in
{\mathbb F}_{q^{dk}}$ for all $1\leq i\leq n$ and $\lambda \in
{\mathbb F}_{q^k}$. Then for every positive integer $k$, we have
the estimate
$$|N_d(k) -(\frac{(q^{kd}-1)^n}{q^{k(d-1)}}
+\frac{1}{2}(1+(-1)^d)q^{k(\frac{d(n-1)}{2}+1)})| \leq
(D+2)q^{k(\frac{d(n-1)+1}{2})},$$ where $D$ is the total degree of
the rational function $P_d(T)$.

\end{cor}

Since the first Hodge number $h^{0, n-1}(X_{\lambda})=1$, the zeta
function of each fibre $X_{\lambda}$ has at most one non-trivial
$p$-adic unit root. One deduces the $p$-adic continuity result: If
$nm+1\leq d_1\leq d_2$ are positive integers such that
$$d_1\equiv
d_2 ~({\rm mod}~ (p-1)p^m),$$ then
$$Z_{d_1}({\mathbb A}^1, X_{\lambda}) \equiv Z_{d_2}({\mathbb A}^1, X_{\lambda}) ~({\rm
mod}~ p^{m+1}).$$ For a $p$-adic integer $s\in \Zp$ and a residue
class $r\in \Z/(p-1)\Z$, let $\{ d_i\}_{i=1}^{\infty}$ be a
sequence of positive integers in the residue class $r$ mod$(p-1)$,
going to infinity as complex numbers but approaching to $s$ as
$p$-adic numbers, then the limit
$$\zeta_{r,s}({\mathbb A}^1, X_{\lambda}) =\lim_{i\rightarrow \infty}
Z_{d_i}({\mathbb A}^1, X_{\lambda}) \in 1+T\Zp[[T]]$$ exists as a
formal $p$-adic power series. This limit depends only on $s$ and
$r$, not on the particular chosen sequence $\{
d_i\}_{i=1}^{\infty}$. The limit $\zeta_{r,s}({\mathbb A}^1,
X_{\lambda})$ is precisely Dwork's unit root zeta function
attached to the family $X_{\lambda}$. It is a $p$-adic meromorphic
function in $T$ for every $s\in \Zp$ and $r\in \Z/(p-1)\Z$, as
conjectured by Dwork \cite{Dw} and proven by Wan \cite{W2}.
It should be viewed as a two variable $p$-adic zeta function
in $(s, T)$.
Similarly, combining $p$-adic methods in \cite{W2} and $\ell$-adic
methods, we show the limit
$${\mathcal P}_{r,s,b}(T) =\lim_{i\rightarrow \infty}
P_{d_i, b}(T) \in 1+T\Zp[[T]]$$ exists and is in fact a $p$-adic
entire function for each $r$, $s$ and $b$. Taking the limit of the
previous theorem, we obtain the following result for Dwork's unit
root zeta function.

\begin{thm}Assume that $p$ does not divide $n+1$. Then, Dwork's
unit root zeta function is given by
$$\zeta_{r,s}({\mathbb A}^1, X_{\lambda})^{(-1)^{n-1}} =
\frac{1-T}{(1-qT)^{n+1}}\prod_{b=0}^n
{\mathcal P}_{r,s,b}(T)^{(-1)^{b-1}(b-1)}.$$
\end{thm}

In the elliptic family case, the $p$-adic entire function
${\mathcal P}_{r,s,0}(T)$ ($b=0$) is the characteristic power series of the
$U_p$-operator acting on the space of overconvergent $p$-adic
modular forms. It would be interesting but apparently difficult to
get precise information on the slopes for the entire function
${\mathcal P}_{r,s,b}(T)$. Buzzard \cite{Bu} had a complicated but explicit
conjecture in some elliptic modular cases.

We now briefly explain the ideas in proving the above theorems.
For a prime $\ell\not=p$, let ${\mathcal H}^j(K)$ denote the relative
$\ell$-adic cohomology with compact support of the family
$X_{\lambda}$. Then, $Z_d({\mathbb A}^1, X_{\lambda})$ can be
expressed in terms of the L-function over ${\mathbb A}^1$ of the
$d$-th Adams operation of the sheaf ${\mathcal H}^j(K)$:
$$Z_d({\mathbb A}^1, X_{\lambda}) = \prod_{j=0}^{2(n-1)} L({\mathbb A}^1,
[{\mathcal H}^j(K)]^d)^{(-1)^j}.$$ It is thus a rational function in
$T$ for each positive integer $d$.

For a prime number $\ell$ which may be equal to $p$,  let ${\mathcal
F}_{\ell}$ be the non-trivial part of the relative $\ell$-adic
cohomology with compact support of the family $X_{\lambda}$
parameterized by $\lambda \in {\mathbb A}^1$. If $\ell \not=p$,
then ${\mathcal F}_{\ell}$ is the non-trivial part of the middle
dimensional relative cohomology ${\mathcal H}^{n-1}(K)$ and the
generic rank of ${\mathcal F}_{\ell}$ is $n$. If $\ell=p$, then the
generic rank of ${\mathcal F}_p$ is $1$ as the first Hodge number
$h^{0, n-1}(X_{\lambda})=1$ and the family $X_{\lambda}$ is
generically ordinary \cite{W31}. The $d$-th moment zeta function
is then given up to trivial factors, by the $d$-th moment
L-function:
$$Z_d({\mathbb A}^1, X_{\lambda}) \sim L({\mathbb A}^1, [{\mathcal
F}_{\ell}]^d)^{(-1)^{n-1}},$$where $[{\mathcal F}_{\ell}]^d$ denotes
the $d$-th Adams operation of the sheaf ${\mathcal F}_{\ell}$ on
${\mathbb A}^1$. Similarly, the unit root zeta function
$\zeta_{r,s}({\mathbb A}^1, X_{\lambda})$ is given up to trivial
factors by the unit root L-function:
$$\zeta_{r,s}({\mathbb A}^1, X_{\lambda}) \sim L({\mathbb A}^1,
\omega({\mathcal F}_p)^r \otimes ({\mathcal F}_p\otimes \omega({\mathcal
F}_p)^{-1})^s)^{(-1)^{n-1}},$$ where $\omega({\mathcal F}_p)$ denotes
the Teichm\"uller lifting of the reduction ${\mathcal F}_p\otimes
{\mathbb F}_p$.

Fix a prime number $\ell \not=p$, let $\mathcal F$ denote the
$\ell$-adic sheaf ${\mathcal F}_{\ell}$. For non-negative integers $a$
and $b$, let
$${\mathcal G}_{a,b}:=\Sym^a{\mathcal F}\otimes\wedge^b{\mathcal F},$$
which is an $\ell$-adic sheaf on ${\mathbb A}^1$, vanishing if
$b>n$. Thus, we shall assume that $0\leq b\leq n$ from now on. The
generic rank of ${\mathcal G}_{a,b}$ is ${{n+a-1}\choose{a}}{n\choose
b}$, which goes to infinity as $a$ goes to infinity. The $d$-th
moment L-function is then given \cite{W20} by the formula
$$L({\mathbb A}^1, [{\mathcal F}_{\ell}]^d) = \prod_{b=0}^n L({\mathbb A}^1, {\mathcal G}_{d-b, b})^{(-1)^{b-1}(b-1)}.$$
Thus, to a large extent, the moment zeta functions are reduced to
the study of the L-function $L({\mathbb A}^1, {\mathcal G}_{a, b})$ of
the sheaf ${\mathcal G}_{a,b}$ for all non-negative integers $a$ and
$b$. To understand the purity decomposition and the trivial
factors of this last L-function, the key is to determine the local
and global monodromy of the sheaf $\mathcal F$. This is accomplished
in Section 2. As a consequence, we obtain

\begin{thm} Assume that $p$ does not divide $n+1$. Let $a$ and $b$
be non-negative integers with $0\leq b\leq n$. Then, we have the
formula
$$
L({\mathbb A}^1, {\mathcal
G}_{a,b})=\frac{P_{a,b}(T)\prod_{k=0}^{[(a+b)(n-1)/2]}(1-q^kT)^{\alpha_{a,b}(k)}}
{(1-q^{(a+b)(n-1)/2}T)^{\delta_{a,b}}(1-q^{(a+b)(n-1)/2+1}T)^{\delta_{a,b}}},
$$
where $P_{a,b}(T)\in 1+T\Z[T]$ is a polynomial whose degree is
explicitly given, $P_{a,b}(T)$ is pure of weight $(a+b)(n-1)+1$,
$\delta_{a,b}=0$ or $1$ is explicitly given by Proposition
\ref{dimension}, $\alpha_{a,b}(k)$ is the coefficient of $x^kz^b$
in the power series
$$
\{\frac{(1-x^n)\cdots(1-x^{a+n-1})}{(1-x^2)\cdots(1-x^a)}\}(1+z)(1+xz)\cdots(1+x^{n-1}z),
$$
where the quantity in the bracket is understood to be $1-x^n$ if
$a=1$, and $1-x$ if $a=0$.
\end{thm}

Let now $s\in \Zp$ be a $p$-adic integer and let $r$ be a residue
class modulo $(p-1)$. Choose a sequence of positive integers $\{
d_i\}_{i=1}^{\infty}$ in the residue class $r$ modulo $(p-1)$,
going to infinity as complex numbers but approaching to $s$ as
$p$-adic integers. For integers $0\leq b\leq n$, we define
$${\mathcal L}_{r,s,b}({\mathbb A}^1, T) =\lim_{i\rightarrow \infty}L({\mathbb A}^1, {\mathcal G}_{d_i-b,b})\in 1+\Zp[[T]].$$ This limit exists as a formal $p$-adic power
series. It depends only on $r, s$ and $b$, not on the choice of
the sequence $\{ d_i\}_{i=1}^{\infty}$. It follows from the
general result in \cite{W2} that ${\mathcal L}_{r,s,b}({\mathbb A}^1, T)$ is
a $p$-adic meromorphic function in $T$. The formula
$$L({\mathbb A}^1, \omega({\mathcal F}_p)^r \otimes ({\mathcal F}_p\otimes \omega({\mathcal
F}_p)^{-1})^s)=\prod_{b=0}^n {\mathcal L}_{r,s,b}({\mathbb A}^1,
T)^{(-1)^{b-1}(b-1)}$$shows that the unit root L-function on the left
side is also $p$-adic meromorphic in $T$. It study is reduced, to
a large extent, to the study of the L-functions
${\mathcal L}_{r,s,b}({\mathbb A}^1, T)$ for all $r, s$ and $b$.

Combining the above theorem together with the $p$-adic limiting
argument in \cite{W2}, we obtain the following more precise
result. The proof is similar to the one given in \cite{fw2} for
the Kloosterman family. The key point is that the number of
$p$-adic zeros of the polynomial $P_{a,b}(T)$ in any fixed
$p$-adic disc $|T|_p <M$ ($M$ finite) is uniformly bounded for all
$a$ and $b$. This fact holds only for those motive ${\mathcal F}$
whose first Hodge number is $1$, which is the case for Calabi-Yau
hypersurfaces.

\begin{thm} Assume that $p$ does not divide $n+1$. Let $a$ and $b$
be non-negative integers with $0\leq b\leq n$. Then, for each
$s\in \Zp$ and each residue class $r\in \Z/(p-1)\Z$, we have the
factorization
$${\mathcal L}_{r,s,b}({\mathbb A}^1, T)= {\mathcal P}_{r,s,b}(T)\prod_{k=0}^{\infty}(1-q^kT)^{\beta_b(k)},$$ where
${\mathcal P}_{r,s,b}(T)\in 1+T\Zp[[T]]$ is a $p$-adic entire function and
$\beta_b(k)$ is the coefficient of $x^kz^b$ in the power series
$$
\frac{(1+z)(1+xz)\cdots(1+x^{n-1}z)}{(1-x^2)(1-x^3)\cdots
(1-x^{n-1})}.
$$
In particular, ${\mathcal L}_{r,s,b}({\mathbb A}^1, T)$ is a $p$-adic entire
function with a zero at $T=q^{-k}$ of multiplicity at least
$\beta_b(k)$ for each non-negative integer $k$.
\end{thm}

It would be interesting to determine the slopes of the polynomials
$P_{a,b}(T)$ and the entire functions ${\mathcal P}_{r,s,b}(T)$. This seems
to be quite difficult in general. The simplest case $n=2$ (the
elliptic family case) has been studied extensively in connection
to slopes of modular forms, over-convergent $p$-adic modular forms
\cite{W1}, the eigencurve \cite{CM} and the Gouvea-Mazur
conjectures. The first step may be to get a good explicit lower
bound for the $p$-adic Newton polygon of $P_{a,b}(T)$ and
${\mathcal P}_{r,s,b}(T)$. Such a good lower bound is already quite
non-trivial to obtain.

Taking partial derivative with respect to $z$ in the generating
function for $\beta_b(k)$ and then setting $z=-1$, we deduce
$$\sum_{k=0}^{\infty}(\sum_{b=0}^n (-1)^{b-1}b\beta_b(k))x^k = 1-x.$$ This together with the previous theorem implies

\begin{cor}
Assume that $p$ does not divide $n+1$. Then, the unit root
L-function is given by
$$L({\mathbb A}^1, \omega({\mathcal F}_p)^r \otimes ({\mathcal F}_p\otimes \omega({\mathcal
F}_p)^{-1})^s)={\frac{1-T}{1-qT}}\prod_{b=0}^n
{\mathcal P}_{r,s,b}(T)^{(-1)^{b-1}(b-1)}.$$
\end{cor}

It would be of great interest to understand the cancellation
nature in the above alternating product of $p$-adic entire
functions.

The paper is organized as follows. In Section $2$, we determine
both the local monodromy and the global monodromy of the sheaf
${\mathcal F}$. These results are then used in Section $3$ to
calculate the L-function of the sheaf ${\mathcal G}_{a,b}$ and its
local factors at bad points. In Section $4$, we treat the
degenerate case when $p$ divides $n+1$.

{\bf Acknowledgements}. We thank L. Fu and N. Katz for helpful comments
and for providing several relevant references. We were informed by Katz that
many of the results in Section $2$ were proved indepedently by him in his forthcoming
paper \cite{k2} on the Dwork family. The second author was partially supported by NSF. The first author was partially supported by MTM2004-07203-C02-01 and FEDER.

\section{The monodromy via Fourier transform.}

Let $k={\mathbb F}_q$ be a finite field of characteristic
$p$, $n\geq 2$ an integer, $X\subset\AAA^{n+1}_k$ the hypersurface
defined by $x_1\cdots x_{n+1}=1$, and $\sigma:X\to\AAA^1_k$ the
restriction of the sum map $(x_1,\ldots,x_{n+1})\to
x_1+\ldots+x_{n+1}$ to $X$. Fix a prime $\ell\neq p$. We want to
study the local monodromy of the non-trivial part of the object
$K:=\R\sigma_!\QQ\in\Dbc$, which parameterizes the cohomology of the family described in the introduction. The main results are summarized in the
following theorem:

\begin{thm}\label{main} The cohomology sheaves ${\mathcal H}^j(K)=\R^j\sigma_!\QQ$ vanish
for $j<n-1$ and $j>2n-2$. We have isomorphisms
$$
{\mathcal H}^j(K)\cong\QQ^{n\choose{j-n+2}}(n-1-j)
$$
for $n\leq j\leq 2n-2$, and an exact sequence
$$
0\to\QQ^n\to{\mathcal H}^{n-1}(K)\to{\mathcal F}\to 0
$$
where $\mathcal F$ is the extension by direct image of a geometrically irreducible
smooth sheaf on the dense open set $U=\AAA^1_k-\{(n+1)\zeta:\zeta^{n+1}=1\}$,
of rank $n$ and punctually pure of weight $n-1$. It is endowed with a
non-degenerate pairing $\Phi:{\mathcal F}\times{\mathcal F}\to\QQ(1-n)$, which is symmetric
if $n$ is odd and skew-symmetric if $n$ is even. As a representation of the inertia
group at infinity, $\mathcal F$ is unipotent with a single Jordan block.

If $p$ does not divide $n+1$, $\mathcal F$ is everywhere tamely
ramified. The inertia group at each of the $n+1$ singular points
$x=(n+1)\zeta$ acts on ${\mathcal F}_{\bar\eta}$ with invariant subspace
of codimension $1$. On the quotient ${\mathcal F}_{\bar\eta}/{\mathcal
F}_{\bar\eta}^{I_x}$, $I_x$ acts trivially if $n$ is even, and
through its unique character of order $2$ if $n$ is odd.

If $p$ divides $n+1$, let $n+1=p^am$, with $m$ prime to $p$. Then
$\mathcal F$ is smooth on $\GG_m$, and the inertia group at $0$ acts
with invariant subspace of dimension $m-1$. The action of $I_0$ on
the quotient ${\mathcal F}_{\bar\eta}/{\mathcal F}_{\bar\eta}^{I_0}$ is
totally wild, with a single break $1/(p^a-1)$ with multiplicity
$m(p^a-1)=n-m+1$. In particular, the Swan conductor at $0$ is $m$.

The determinant of $\mathcal F$ is the geometrically constant sheaf $\QQ(-n(n-1)/2)$ if $n$ is even or $p$ divides $n+1$, and the pulled back Kummer sheaf $${\mathcal L}_{\chi(\lambda^{n+1}-(n+1)^{n+1})}(-n(n-1)/2)$$ if $n$ is odd and $(p,n+1)=1$, where $\chi$ is the unique character of order $2$ of the inertia group $I_0$.

The geometric monodromy group of $\mathcal
F$ is given by

$$
\left\{
\begin{array}{ll} Sp(n,\Phi) & \mbox{if $n$ is even} \\
O(n,\Phi) & \mbox{if $n$ is odd and $(p,n+1)=1$}
\\ SO(n,\Phi) & \mbox{if $n$ is odd, $p|n+1$}  
\\ & \mbox{  and $(p,n)\neq (2,5)$ or $(2,7)$}
\\ G_2  \mbox{ in its standard} & 
\\ \mbox{    7-dimensional representation} & \mbox{if $p=2$, $n=7$}
\\ SL(2) \mbox{ in $sym^4$ of its} &
\\ \mbox{    standard representation} & \mbox{if $p=2$, $n=5$}
\end{array}\right.
$$
\end{thm}

\bigskip

We will deduce most of the properties of the object $K$ from the
properties of its Fourier transform $L\in\Dbc$ with respect to a
fixed non-trivial additive character $\psi:k\to{\mathbb
C}^\star\isomto\QQ^\star$. The Fourier transform $L$ is closely
related to the Kloosterman sheaf. This connection of the Dwork family with
Kloosterman sums was first discovered by Katz \cite{ka1} (Section 5.5)
who uses the properties of the family to get information on certain Kloosterman
sums. We will use this connection the other way around and apply
Katz's fundamental results for the Kloosterman sheaf.

Recall (cf. \cite{kl}) that the Fourier transform is defined by
$$
FT_{\psi}(K)=\R\pi_{2!}(\pi_1^\star(K)\otimes\mu^\star{\mathcal L}_\psi)[1]
$$
where $\pi_1,\pi_2:\AAA^2_k\to\AAA^1_k$ are the projections,
$\mu:\AAA^2_k\to\AAA^1_k$ is the product map and ${\mathcal L}_\psi$
is the Artin-Schreier sheaf on $\AAA^1_k$ associated to the
character $\psi$. It is an auto-equivalence of the triangulated
category $\Dbc$, and has the following involution property:
$FT_{\bar\psi}FT_{\psi}(K)=K(-1)$.

One of the main advantages of this equivalence is that, following Laumon (cf. \cite{laumon}), the local properties of the object $K$ can be read from those of its Fourier transform. This is the method that we will use to deduce most of the results about $K$.

Let us first determine what the Fourier transform of $K$ is explicitly. Using proper base change on the cartesian diagram
$$\begin{CD}
X @<{\tilde\pi_1}<< X\times\AAA^1_k \\
@V{\sigma}VV @VV{\tilde\sigma}V\\
\AAA^1_k @<{\pi_1}<< \AAA^2_k
\end{CD}$$
we get
$$
\pi_1^\star(K)=\pi_1^\star(\R\sigma_!\QQ)=\R\tilde\sigma_!\tilde\pi_1^\star\QQ=\R\tilde\sigma_!\QQ
$$
By the projection formula, we have then
$$
L=\R\pi_{2!}((\R\tilde\sigma_!\QQ)\otimes\mu^\star{\mathcal
L}_\psi)[1]=\R\pi_{2!}(\R\tilde\sigma_!(\tilde\sigma^\star\mu^\star{\mathcal
L}_\psi))[1]=\R\tilde\pi_{2!}(\tilde\mu^\star{\mathcal L}_\psi)[1]
$$
where $\tilde\pi_1$ and $\tilde\pi_2$ are the projections of
$X\times\AAA^1_k$ onto its factors and
$\tilde\mu:X\times\AAA^1_k\to\AAA^1_k$ is the map
$((x_1,\ldots,x_{n+1}),t)\mapsto t(x_1+\ldots+x_{n+1})$.

Extend the canonical map $L\to j_\star j^\star L$ to a distinguished triangle \begin{equation}\label{triangle}
 M\to L\to j_\star j^\star L\to.
\end{equation}
in $\Dbc$, where $j:\AAA^1_k-\{0\}\hookrightarrow\AAA^1_k$ is the open
immersion. The object $M$ is punctual supported at $0$, since $L\to j_\star j^\star L$ is an isomorphism away from $0$.

At $0$, the object $L$ is just $\R\Gamma_c(X\otimes \bar k,\QQ)[1]$
by proper base change. Since $X$ is just the product of $n$ copies
of $\GG_m$, we have
$$
L_0=\bigotimes_{i=1}^n \R\Gamma_c(\GG_{m,\bar k},\QQ)[1]
$$
From $\HH^1_c(\GG_{m,\bar k},\QQ)=\QQ$, $\HH^2_c(\GG_{m,\bar
k},\QQ)=\QQ(-1)$ and $\HH^i_c(\GG_{m,\bar k},\QQ)=0$ for $i\neq
1,2$, we conclude
$$
{\mathcal H}^{i-1}(L)_0=\QQ^{n\choose{i-n}}(n-i)
$$
for $n\leq i\leq 2n$, and $0$ otherwise, so we get a quasi-isomorphism
$$
L_0\cong\bigoplus_{i=n}^{2n}\QQ^{n\choose{i-n}}(n-i)[1-i]
$$

Away from $0$, we have $L=\R\tilde\pi_{2!}(\tilde\mu^\star{\mathcal
L}_\psi)[1]$, where we now regard $\tilde\pi_{2!}$ as the projection
$X\times\GG_m\to\GG_m$. Consider the
automorphism $\phi$ of $\AAA^{n+1}_k\times\GG_m$ given by
$\phi((x_1,\ldots,x_{n+1}),t)=((tx_1,\ldots,tx_{n+1}),t)$. The image
of $X\times\GG_m$ under $\phi$ is the variety $Y$ defined by the
equation $x_1\cdots x_{n+1}=t^{n+1}$, and
$\tilde\mu=\tilde\sigma\circ\phi$. Since $\phi$ is an automorphism,
$\phi^\star=\R\phi_\star=\R\phi_!$, and we get
$$
j^\star L=\R\tilde\pi_{2!}(\tilde\mu^\star{\mathcal
L}_\psi)[1]=\R\tilde\pi_{2!}(\phi^\star\tilde\sigma^\star{\mathcal
L}_\psi)[1]=
$$
$$
=\R(\tilde\pi_2\phi)_!(\tilde\sigma^\star{\mathcal
L}_\psi)[1]=\R\tilde\pi_{2!}(\tilde\sigma^\star{\mathcal L}_\psi)[1]
$$
The stalk of $j^\star L$ at a geometric point $t\in\GG_{m,\bar k}$
is then $\R\Gamma_c(\{x_1\cdots x_{n+1}=t^{n+1}\},\LL_{\psi(\sum
x_i)})[1]$. By \cite{deligne}, Th\'eor\`eme 7.4, we deduce that
${\mathcal H}^i(j^\star L)=0$ for $i\neq n-1$, and ${\mathcal
H}^{n-1}(j^\star L)$ is the pull-back by the $(n+1)$-th power map of
the Kloosterman sheaf given in \cite{deligne}, Th\'eor\`eme 7.8 and,
more generally, in \cite{gkm}, 4.1.1. Therefore we have a quasi-isomorphism
$$
j^\star L\cong[n+1]^\star\mathrm{Kl}_{n+1}(\psi)[1-n]
$$
Denote by $\mathcal L$ the sheaf $[n+1]^\star\mathrm{Kl}_{n+1}(\psi)$ on
$\GG_m$. It is geometrically irreducible, because it is already
irreducible as a representation of the inertia group at $0$: by
\cite{deligne}, Th\'eor\`eme 7.8, the action of a topological generator is unipotent with a
single Jordan block. In particular, the invariant subspace for the inertia action at $0$ has dimension $1$, so the stalk of $j_\star j^\star L$ at $0$ is quasi-isomorphic to $\QQ[1-n]$.

Taking stalks at $0$ in the distinguished triangle (\ref{triangle}) we get
$$
M_0\to \bigoplus_{i=n}^{2n}\QQ^{n\choose{i-n}}(n-i)[1-i] \to \QQ[1-n] \to.
$$
and consequently a quasi-isomorphism
$$
M_0\simeq \bigoplus_{i=n+1}^{2n}\QQ^{n\choose{i-n}}(n-i)[1-i]=\bigoplus_{i=1}^{n}\QQ^{n\choose{i}}(-i)[1-n-i]
$$
Then since $M$ is punctual supported at $0$, the distinguished triangle (\ref{triangle}) reads
$$
\bigoplus_{i=1}^{n}\QQ^{n\choose{i}}(-i)[1-n-i]_0\to L\to j_\star{\mathcal L}[1-n]\to.
$$

Taking Fourier transform with respect to the complex conjugate
character $\bar\psi$ and using the facts that
$FT_{\bar\psi}FT_\psi(K)=K(-1)$ and that the Fourier transform of the punctual sheaf $(\QQ)_0$ is the shifted constant sheaf $\QQ[1]$, we get the distinguished triangle
$$
\bigoplus_{i=1}^{n}\QQ^{n\choose{i}}(-i)[2-n-i]\to K(-1)\to FT_{\bar\psi}(j_\star{\mathcal L})[1-n]\to.
$$
Since $\mathcal L$ is a geometrically irreducible sheaf of rank $\geq 2$, its direct image $j_\star{\mathcal L}$ is a Fourier sheaf in the sense of \cite{gkm}, 8.2 (cf. \cite{gkm}, lemma 8.3.1). Then its Fourier transform is a sheaf of the same kind, by (\cite{gkm}, Theorem 8.2.5). Namely, it is the extension by direct image to $\AAA^1$ of a geometrically irreducible sheaf on a dense open set $U\subset\AAA^1$, and we get a distinguished triangle
$$
\bigoplus_{i=1}^{n}\QQ^{n\choose{i}}(1-i)[2-n-i]\to K\to {\mathcal F}[1-n]\to.
$$
where ${\mathcal F}=FT_{\bar\psi}(j_\star{\mathcal L})(1)$. Taking the associated long exact sequence of cohomology sheaves and using the fact that ${\mathcal F}$ has no punctual sections, we get an exact sequence
$$
0\to \QQ^n\to {\mathcal H}^{n-1}(K)\to {\mathcal F} \to 0
$$
and isomorphisms
$$
{\mathcal H}^j(K)\isomto \QQ^{n\choose{j-n+2}}(n-1-j)\mbox{ for }n\leq
j\leq 2n-2
$$
and
$$
{\mathcal H}^j(K)=0 \mbox{ for }j\not\in\{n-1,\ldots,2n-2\}.
$$

Thus the cohomology of our family has a
``constant part", which has dimension $n\choose{j-n+2}$ and is
pure of weight $2(j-n+1)$ on degree $j$ for every
$j={n-1},\ldots,2n-2$, and a non-constant geometrically
irreducible part on degree $n-1$ given by the sheaf $\mathcal F$. If $n+1$ is prime to $p$, this sheaf is the pull-back by the $(n+1)$-th power map of a hypergeometric sheaf, as defined by Katz in \cite{esde} (Section 8). Namely, using the same notation as in the reference, it is $[n+1]^\star\mathrm{Hyp}_{(n+1)^{n+1}}(!,\psi,$all nontrivial characters $\chi$ of order dividing $n+1$;$n$ times the trivial character$)$ (cf. \cite{esde}, Theorem 9.3.2). We will not make use of this fact in what follows.

\begin{prop} The sheaf $\mathcal F$ is smooth of rank $n$ and punctually
pure of weight $n-1$ on
$U=\AAA^1_{\bar k}-\{(n+1)\zeta:\zeta^{n+1}=1\}$.\end{prop}

{\bf Proof.} If $n+1$ is prime to $p$, by \cite{fw}, lemma 1.4, the wild inertia group of
$\AAA^1_{\bar k}$ at infinity acts on $\mathcal L$ as
$\bigoplus_{\zeta^{n+1}=1}\LL_{\psi_{(n+1)\zeta}}$, where
$\psi_{(n+1)\zeta}(t)=\psi((n+1)\zeta t)$. By \cite{esde}, Lemma 7.3.9, $\mathcal F$ is
smooth at $t\in\AAA^1_{\bar k}$ if and only if all breaks of ${\mathcal
L}\otimes\LL_{\bar\psi_t}$ at infinity are $\geq 1$. But, as a representation of
$P_\infty$, ${\mathcal
L}\otimes\LL_{\bar\psi_t}=\bigoplus_{\zeta^{n+1}=1}\LL_{\psi_{(n+1)\zeta-t}}$ has all
its breaks equal to $1$ unless $t=(n+1)\zeta$ for some $\zeta\in\mu_{n+1}({\bar
k})$. This proves that $\mathcal F$ is smooth on $U$. If $p$ divides $n+1$, all breaks of
$\mathcal L$ at infinity are $<1$, so $\mathcal F$ is smooth on $U=\GG_{m,{\bar k}}$ by
\cite{gkm}, 8.5.8.

Since $\mathcal L$ is pure of weight $n$, so is its direct image $j_\star{\mathcal L}[0]$ as a derived category object. The Fourier transform
preserves purity and shifts weights by $1$, so ${\mathcal F}(-1)[0]$ is pure of weight $n+1$ as a derived category object. In particular, on the open set where $\mathcal F$ is smooth, it is punctually pure of
weight $(n+1)-2=n-1$. To compute the rank, we
use the formula in \cite{esde}, 7.3.9, which gives
$$
\mathrm{rank}({\mathcal F})=\mathrm{drop}_0({\mathcal L})=(n+1)-1=n.
$$\hfill$\Box$

\begin{prop}\label{pairing} There is a non-degenerate pairing $\Phi:{\mathcal
F}_U\times{\mathcal F}_U\to\QQ(1-n)$ which is symmetric for $n$ odd and
skew-symmetric for $n$ even.
\end{prop}

{\bf Proof.} According to \cite{gkm} 4.1.3, the dual of the sheaf
$\mathrm{Kl}_{n+1}(\psi)$ on $\GG_{m,k}$ is
$\mathrm{Kl}_{n+1}(\bar\psi)(n-1)$. Therefore, the dual of the
object $j_\star{\mathcal L}[0]\in\Dbc$ is
$j_\star\bar{\mathcal L}[0](n-1)$, where $\bar{\mathcal
L}=[n+1]^\star\mathrm{Kl}_{n+1}(\bar\psi)$.

By \cite{kl}, Th\'eor\`eme 2.1.5, the dual of the Fourier
transform with respect to $\psi$ of an object is the Fourier
transform with respect to $\bar\psi$ of the dual object.
Therefore, the dual of $FT_{\bar\psi}(j_\star{\mathcal L}[0])={\mathcal
F}[0]$ is $FT_\psi(j_\star\bar{\mathcal L}[0](n-1))={\mathcal F}[0](n-1)$.
In particular, we have a non-degenerate pairing on the open set $U$ where $\mathcal F$ is smooth: ${\mathcal
F}_U\times{\mathcal F}_U\to\QQ(1-n)$. Since ${\mathcal F}_U$ is
irreducible, the pairing is unique up to a scalar and either symmetric of skew-symmetric.
The actual sign is given by the usual cup product sign, since
$\mathcal F$ is a subsheaf of $\R^{n-1}\sigma_!(\QQ)$.\hfill$\Box$

\begin{prop}\label{unipotent} The sheaf $\mathcal F$ is tamely ramified at infinity. The tame
inertia group at infinity $I^{tame}_\infty$ acts unipotently on
${\mathcal F}_{\bar\eta}$ with a single Jordan block.\end{prop}

{\bf Proof.} Since $\mathcal L$ is tamely ramified at $0$
and the inertia group acts unipotently with a single Jordan block,
the same is true for $\mathcal F$ at $\infty$ by \cite{esde}, Theorem
7.5.4.\hfill$\Box$.

\begin{prop} \label{finitepoints} Suppose that $n+1$ is prime to $p$. Then $\mathcal F$ is
everywhere tamely ramified, and for every $(n+1)$-th root of unity
$\zeta$ in $\bar k$, the action of the inertia group at $(n+1)\zeta$
on ${\mathcal F}_{\bar\eta}$ has invariant subspace of codimension $1$.
\end{prop}

{\bf Proof.} Let $\zeta$ be a $n+1$-th root of unity in $\bar k$.
Then $\zeta:(x_1,\ldots,x_{n+1})\to(\zeta x_1,\ldots,\zeta x_{n+1})$
is an automorphism of $X$. Therefore,
$K:=\R\sigma_!\QQ=\R(\sigma\circ\zeta)_!\QQ=\R(\tilde\zeta\circ\sigma)_!\QQ=[\tilde\zeta]_\star\R\sigma_!\QQ=[\tilde\zeta]_\star
K$ where $[\tilde\zeta]:\AAA^1_k\to\AAA^1_k$ is multiplication by
$\zeta$. So the sheaf $\mathcal F$ is invariant under multiplication by
$(n+1)$-th roots of unity on $\AAA^1_k$. In particular, the local
monodromies at $(n+1)\zeta$ are isomorphic for all
$\zeta\in\mu_{n+1}(\bar k)$.

By the Euler-Poincar\'e formula,
$$
\chi_c({\mathcal F})=\mathrm{rank}({\mathcal
F})-\sum_{t\in(n+1)\mu_{n+1}(\bar k)}(\mathrm{drop}_t{\mathcal
F}+\mathrm{swan}_t{\mathcal F})
$$
since $\mathcal F$ is tamely ramified at infinity and smooth on
$\AAA^1_k-(n+1)\mu_{n+1}(\bar k)$. We can compute this Euler characteristic
directly:
$$
\chi_c(K)=\chi_c(\R\sigma_!{\QQ})=\chi_c(X,\QQ)=0
$$
since $X$ is a product of copies of $\GG_m$. Therefore
$$
0=\chi_c(K)=\sum_{j=n-1}^{2n-2}(-1)^j\chi_c({\mathcal H}^j(K))$$
$$
=(-1)^{n-1}\chi_c({\mathcal
F})+(-1)^{n-1}n+\sum_{j=n}^{2n-2}(-1)^j{n\choose{j-n+2}}$$
$$
=(-1)^{n-1}\chi_c({\mathcal F})+\sum_{j=1}^n (-1)^{j+n}{n\choose
j}=(-1)^{n-1}\chi_c({\mathcal F})-(-1)^n,
$$
so $\chi_c({\mathcal F})=-1$. We conclude that
$$
\sum_{t\in(n+1)\mu_{n+1}(\bar k)}(\mathrm{drop}_t{\mathcal
F}+\mathrm{swan}_t{\mathcal F})=n+1
$$
and therefore the only possibility is $\mathrm{drop}_t{\mathcal F}=1$
and $\mathrm{swan}_t{\mathcal F}=0$ for every $t\in(n+1)\mu_{n+1}({\bar
k})$. In particular, $\mathcal F$ is everywhere tamely ramified.\hfill
$\Box$

\begin{prop} \label{root} Suppose that $n+1$ is prime to $p$, and let $t\in(n+1)\mu_{n+1}({\bar
k})$. If $n$ is even, the inertia group $I_t$ acts trivially on the
one-dimensional space ${\mathcal F}_{\bar\eta}/{\mathcal
F}_{\bar\eta}^{I_t}$. That is, the action of $I_t$ on ${\mathcal
F}_{\bar\eta}$ is unipotent with a Jordan block of size $2$ and all
other blocks of size $1$. If $t\in{\mathbb F}_q$, the action of a
geometric Frobenius element at $t$ on ${\mathcal F}_{\bar\eta}^{I_t}$
has one of $\pm q^{(n-2)/2}$ as an eigenvalue, and all other
eigenvalues of absolute value $q^{(n-1)/2}$. 

If $n$ is odd, $I_t$
acts on the one-dimensional space ${\mathcal F}_{\bar\eta}/{\mathcal
F}_{\bar\eta}^{I_t}$ via its unique character of order $2$. In
particular, the action of $I_t$ on ${\mathcal F}_{\bar\eta}$ is
semisimple. If $t\in{\mathbb F}_q$, the action of a geometric
Frobenius element at $t$ on ${\mathcal F}_{\bar\eta}^{I_t}$ has all
eigenvalues of absolute value $q^{(n-1)/2}$. \end{prop}

{\bf Proof.} This can be proven using the Picard-Lefschetz formulas
(cf. \cite{sga7ii}, expos\'e XV), since the fibres of
$\sigma:X\to\AAA^1$ have only isolated ordinary quadratic
singularities. Alternatively, one may use the explicit description
of the monodromy at infinity of the Kloosterman sheaf and Laumon's
local Fourier transform theory.

According to \cite{fw}, Theorem 1.1, the action of the inertia group
at infinity on $\mathrm{Kl}_{n+1}(\psi)$ is given by
$[n+1]_\star{\mathcal L}_{\psi_{n+1}}$ if $n$ is even and
$[n+1]_\star({\mathcal L}_{\psi_{n+1}}\otimes{\mathcal L}_{\chi_2})$ if $n$
is odd. Therefore, the action on
$[n+1]^\star\mathrm{Kl}_{n+1}(\psi)$ is given by
$\oplus_{\zeta^{n+1}=1}{\mathcal L}_{\psi_{(n+1)\zeta}}$ if $n$ is even
and $\oplus_{\zeta^{n+1}=1}{\mathcal L}_{\psi_{(n+1)\zeta}}\otimes{\mathcal
L}_{\chi_2}$ if $n$ is odd (cf. \cite{fw}, Lemma 1.4). We conclude
by \cite{esde}, 7.4.1 and 7.5.4.

In particular, the Frobenius eigenvalues of ${\mathcal
F}_{\bar\eta}^{I_t}$ all have weight $n-1$ if $n$ is odd by
\cite{gkm}, 7.0.8. If $n$ is even, there are $n-2$ eigenvalues of
weight $n-1$ and one of weight $n-2$. Since the local $L$-function
has integral coefficients, the non-real eigenvalues must appear in
complex conjugate pairs, and therefore the one with weight $n-2$
must be real, necessarily $\pm q^{(n-2)/2}$.\hfill$\Box$

\begin{prop}\label{div} Suppose that $p$ divides $n+1$, and write $n+1=p^am$
with $(p,m)=1$. Then the inertia group at $0$ acts with invariant
subspace of dimension $m-1$, and its action on the quotient ${\mathcal
F}_{\bar\eta}/{\mathcal F}_{\bar\eta}^{I_0}$ is totally wild, with a
single break $1/(p^a-1)$ with multiplicity
$m(p^a-1)=n-m+1$.\end{prop}

{\bf Proof.} In this case ${\mathcal L}=j_\star [n+1]^\star
\mathrm{Kl}_{n+1}(\psi) = j_\star [m]^\star
[p^a]^\star\mathrm{Kl}_{n+1}(\psi)=j_\star
[m]^\star\mathrm{Kl}_{n+1}(\psi')$, where $\psi'$ is the additive character given by $\psi'(t)=\psi(t^{p^a})$. We deduce by \cite{gkm}, 1.13.1
that $\mathcal L$ is totally wild at $\infty$ with a single break
$m/(n+1)<1$ with multiplicity $n+1$. Therefore, by \cite{esde}
7.5.4, we conclude that $\mathcal F$ has break $m/(n-m+1)$ at $0$ with
multiplicity $n-m+1$. In particular, the Swan conductor at $0$ is
$m$.

It remains to compute the tame part of the monodromy at $0$. By the
Euler-Poincar\'e formula,
$$
-1=\chi_c({\mathcal F})=\dim{\mathcal
F}_{\bar\eta}^{I_0}-\mathrm{swan}_0{\mathcal F}=\dim{\mathcal
F}_{\bar\eta}^{I_0}-m
$$
so $\dim{\mathcal F}_{\bar\eta}^{I_0}=m-1$ which is precisely the
codimension of the wild part. Therefore, the inertia group at $0$
has dimension $m-1$ invariant subspace, and the action in ${\mathcal
F}_{\bar\eta}/{\mathcal F}_{\bar\eta}^{I_0}$ is totally wild, with a
single break $m/(n-m+1)=1/(p^a-1)$ with multiplicity $n-m+1$.\hfill$\Box$

\begin{prop}\label{lfunction} The $L$-function of $\mathcal F$ on $\AAA^1_k$ is given by
$$
L(\AAA^1,{\mathcal F},T)=1-T.
$$
The eigenvalues of a geometric Frobenius element
$F_\infty$ at infinity acting on $\mathcal F$ are $1,q,\ldots,q^{n-1}$.\end{prop}

{\bf Proof.} By Theorem \ref{main} we have
$$
L(\AAA^1,K,T)=\prod_{j=n-1}^{2n-2} L(\AAA^1,{\mathcal H}^j(K),T)^{(-1)^j}= $$ $$
=\prod_{j=n-1}^{2n-2} (1-q^{j+2-n}T)^{(-1)^{j+1}{n\choose{j+2-n}}} \cdot L(\AAA^1,{\mathcal
F},T)^{(-1)^{n-1}}
 = $$ $$ =\prod_{j=1}^{n}
(1-q^{j}T)^{(-1)^{j+n-1}{n\choose j}} \cdot L(\AAA^1,{\mathcal F},T)^{(-1)^{n-1}}
$$
On the other hand, we have
$L(\AAA^1,K,T)=L(\AAA^1,R\sigma_!\QQ,T)=Z(X,T)$. Since $X$ is a
product of $n$ copies of the torus $\GG_{m}$, we get
$$
L(\AAA^1,K,T)=\prod_{j=0}^n  (1-q^{j}T)^{(-1)^{j+n-1}{n\choose j}}. $$ Comparing both
expressions, we conclude that $L(\AAA^1,{\mathcal F},T)=1-T$.

Let $j:U\to\PP^1$ be the inclusion. Since $\mathcal F$ is irreducible and not geometrically constant, $H^0(\PP^1,j_\star{\mathcal F})=H^2(\PP^1,j_\star{\mathcal F})=0$. On the other hand, the Euler-Poincar\'e formula gives $\chi(\PP^1,j_\star{\mathcal F})=n+1-\sum_{\zeta^{n+1}=1}1=0$ if $n+1$ is prime to $p$, and $\chi(\PP^1,j_\star{\mathcal F})=1+\dim{\mathcal F}^{I_0}-\mathrm{Sw}_0{\mathcal F}=1+(m-1)-m=0$ if $p$ divides $n+1$, so in either case $H^1(\PP^1,j_\star{\mathcal F})=0$ too. Therefore, the $L$-function of $j_\star{\mathcal F}$ on $\PP^1$ is trivial, so
$$
L(\AAA^1,{\mathcal F},T)=L(\PP^1,j_\star{\mathcal F},T)\det(1-T F_\infty|{\mathcal F}^{I_\infty})= \det(1-T F_\infty|{\mathcal F}^{I_\infty}).
$$
In particular, the action
of $D_\infty/I_\infty$ on the one dimensional space ${\mathcal F}^{I_\infty}$ is trivial, and
the eigenvalues of a geometric Frobenius element acting on ${\mathcal F}$ are
$1,q,\ldots,q^{n-1}$ by \cite{gkm}, 7.0.7.\hfill$\Box$

\begin{prop}\label{determinant} If $n$ is even or $p$ divides $n+1$, the determinant of $\mathcal F$ is the geometrically constant sheaf $\QQ(-n(n-1)/2)$. If $n$ is odd and $(p,n+1)=1$,  $$\det({\mathcal F})={\mathcal L}_{\chi(t^{n+1}-(n+1)^{n+1})}(-n(n-1)/2)$$ where $\chi$ is the unique character of order $2$ of the inertia group of ${\mathbb A}^1$ at $0$ and ${\mathcal L}_{\chi(t^{n+1}-(n+1)^{n+1})}$ is the pullback of the extension by zero to $\AAA^1$ of the corresponding Kummer sheaf on ${\mathbb G}_m$ under the map $t\mapsto t^{n+1}-(n+1)^{n+1}$.\end{prop}

{\bf Proof.} If $n$ is even, the Tate-twisted sheaf ${\mathcal F}((n-1)/2)$ is symplectically self-dual, so its determinant is trivial. Since $\det ({\mathcal F}((n-1)/2))=(\det {\mathcal F}) (n(n-1)/2)$, we conclude that $\det {\mathcal F} \cong\QQ(-n(n-1)/2)$.

If $p$ divides $n+1$, let $n+1=p^am$ as in Proposition \ref{div}. If $\zeta$ is a primitive $m$-th root of unity, exactly as in the proof of Proposition \ref{finitepoints} we get an isomorphism ${\mathcal F}\cong [t\to \zeta t]^\star{\mathcal F}$. In particular, there is a sheaf ${\mathcal G}$ on $\GG_m$ such that ${\mathcal F}_{|\GG_m}=[m]^\star{\mathcal G}$, where $[m]:\GG_m\to\GG_m$ is the $m$-th power map. By \cite{gkm}, 1.13.1 and Proposition \ref{div}, as a representation of the wild inertia group at $0$ the sheaf ${\mathcal G}$ has a single positive break $1/m(p^a-1)=1/(n+1-m)$ with multiplicity $n+1-m=m(p^a-1)>1$, and Swan conductor $1$. At infinity the inertia group acts quasi-unipotently with a single Jordan block, and after tensoring with a suitable Kummer sheaf we can assume that the action is unipotent. Then $\det {\mathcal G}$ is smooth of rank 1 on $\GG_m$, unramified at infinity and its break at $0$ is $\leq 1/(n+1-m)<1$. Since this break (which is the Swan conductor of $\det{\mathcal G}$ at 0) is an integer, it has to be zero. Thus $\det{\mathcal G}$ is tamely ramified at zero, and therefore geometrically trivial, and the same is true for $\det{\mathcal F}=[m]^\star\det{\mathcal G}$.

So there is some $\ell$-adic unit $\alpha$ such that $\det{\mathcal F}\cong\alpha^{\deg}$, where $\alpha^{\deg}$ is the pullback to $\pi_1(\GG_{m,k})$ of the character of $\pi_1(\GG_{m,k})/\pi_1(\GG_{m,\bar k})\cong \mathrm{Gal}(\bar k/k)$ that maps the canonical generator $F$ to $\alpha$. To find the value of $\alpha$ we need to compute the determinant of the action of an element of degree $1$ of $\pi_1(\GG_{m,k})$ on $\det{\mathcal F}$. But from Proposition \ref{lfunction} we know that the action of the geometric Frobenius element at infinity (which has degree $1$) on $\mathcal F$ has eigenvalues $1,q,\ldots,q^{n-1}$. Therefore, $\alpha=q^{1+2+\cdots+(n-1)}=q^{n(n-1)/2}$, so $\det{\mathcal F}\cong(q^{n(n-1)/2})^{\deg}=\QQ(-n(n-1)/2)$.

If $n$ is odd and $(p,n+1)=1$, from Propositions \ref{unipotent} and \ref{root} we know that $\det{\mathcal F}$ is smooth on $U=\AAA^1_{\bar k}-\{(n+1)\zeta:\zeta^{n+1}=1\}$, unramified at infinity and tamely ramified at the $n+1$ singular points $(n+1)\zeta$, with the inertia groups acting via their character $\chi$ of order two. Therefore, $(\det{\mathcal F})\otimes \check{\mathcal L}_{\chi(t^{n+1}-(n+1)^{n+1})}$ is everywhere unramified, and thus geometrically trivial. So there is some $\ell$-adic unit $\alpha$ such that $\det{\mathcal F}\cong\alpha^{\deg}\otimes{\mathcal L}_{\chi(t^{n+1}-(n+1)^{n+1})}$. To find the exact value of $\alpha$, we again evaluate the determinant at $t=\infty$ to be $q^{n(n-1)/2}$ using Proposition \ref{lfunction}. On the other hand, using that ${\mathcal L}_{\chi(t^{n+1}-(n+1)^{n+1})}={\mathcal L}_{\chi(t^{n+1})}\otimes{\mathcal L}_{\chi(1+\frac{(n+1)^{n+1}}{t^{n+1}})}={\mathcal L}_{\chi(1+\frac{(n+1)^{n+1}}{t^{n+1}})}$ (since $\chi$ has order $2$ and $n+1$ is even), we conclude that the Frobenius element at infinity acts trivially on ${\mathcal L}_{\chi(t^{n+1}-(n+1)^{n+1})}$, and therefore $\alpha=q^{n(n-1)/2}$ and $\det({\mathcal F})={\mathcal L}_{\chi(t^{n+1}-(n+1)^{n+1})}(-n(n-1)/2)$.
\hfill$\Box$

\begin{cor} Suppose that $n$ is odd and $(p,n+1)=1$, and let $t\in{\mathbb F}_q$. Then the action of a geometric Frobenius element $F_t$ at $t$ on $\mathcal F$ has $\chi(t^{n+1}-(n+1)^{n+1})q^{(n-1)/2}$ as an eigenvalue (where $\chi:{\mathbb F}_q^\star\to{\mathbb C}^\star$ is the unique character of order $2$) and the remaining eigenvalues appear in complex conjugate pairs.\end{cor}

{\bf Proof.} From the previous theorem we know that the product of the eigenvalues is $\chi(t^{n+1}-(n+1)^{n+1})q^{n(n-1)/2}$. They all have absolute value $q^{(n-1)/2}$ and, given that ${\mathcal F}((n-1)/2)$ is self-dual, they are permuted by the map $z\mapsto q^{n-1}/z$. So the non-real eigenvalues show up in complex conjugate pairs. There are an odd number of real eigenvalues, all of them necessarily equal to $q^{(n-1)/2}$ or $-q^{(n-1)/2}$. Grouping them in pairs of identical eigenvalues, we are left with just one, whose sign must be $\chi(t^{n+1}-(n+1)^{n+1})$ (since the product of the other ones is positive).\hfill$\Box$

\begin{prop} The geometric monodromy group $G$ of $\mathcal F$ is given by
$$
\left\{
\begin{array}{ll} Sp(n,\Phi) & \mbox{if $n$ is even} \\
O(n,\Phi) & \mbox{if $n$ is odd and $(p,n+1)=1$}
\\ SO(n,\Phi) & \mbox{if $n$ is odd, $p|n+1$}
\\ & \mbox{  and $(p,n)\neq (2,5)$ or $(2,7)$}
\\ G_2  \mbox{ in its standard} &
\\ \mbox{    7-dimensional representation} & \mbox{if $p=2$, $n=7$}
\\ SL(2) \mbox{ in $sym^4$ of its} &
\\ \mbox{    standard representation} & \mbox{if $p=2$, $n=5$}
\end{array}\right.
$$
\end{prop}

{\bf Proof.} The connected component $G_0$ of $G$ containing the identity is semisimple by
\cite{weilii}, 1.3.9. Since $G$ contains a unipotent element with a
single Jordan block, its Lie algebra $\mathfrak g$ is simple and
contains a nilpotent element with a single Jordan block and the
representation ${\mathfrak g}\to\mathrm{End}({\mathcal F}_{\bar\eta})$
is faithful and irreducible, by \cite{gkm}, 11.5.2.3. By
\ref{pairing}, we have an a priori inclusion $G\subset Sp(n,\Phi)$
for $n$ even and $G\subset O(n,\Phi)$ for $n$ odd.

Suppose that $n+1$ is prime to $p$. Then $G$ contains
pseudo-reflections (i.e. elements with invariant subspace of
codimension 1). Since any element in $G$ normalizes $\mathfrak g$,
from \cite{esde}, Theorem 1.5 we conclude that ${\mathfrak
g}={\mathfrak{sp}}_n$ if $n$ is even and ${\mathfrak g}={\mathfrak
{so}}_n$ if $n$ is odd. Consequently, $G=Sp(n,\Phi)$ if $n$ is even
and $G=SO(n,\Phi)$ or $O(n,\Phi)$ if $n$ is odd. But the local
monodromies at the points $t\in (n+1)\mu_{n+1}({\bar k})$ contain
elements of determinant $-1$, so $G$ must be the full orthogonal
group.

When $p$ divides $n+1$, we will make use of the classification
theorem in \cite{gkm}, 11.6. According to it, the possibilities for
$\mathfrak g$ are: $\mathfrak{sl}_2$ in the $(n-1)$-th symmetric
power of its standard representation, $\mathfrak{sp}_n$ if $n$ is
even, $\mathfrak{so}_n$ if $n$ is odd and $\mathfrak{g}_2$ in its standard 7-dimensional representation if
$n=7$.

Suppose that $\mathfrak{g}=\mathfrak{sl}_2$, and let $n+1=p^am$ with $m$ prime to $p$.  As in the proof of Proposition \ref{determinant} we find a smooth sheaf $\mathcal G$ on $\GG_m$ such that ${\mathcal F}_{|\GG_m}=[m]^\star{\mathcal G}$. Since the geometric monodromy group of $\mathcal F$ has finite index in that of $\mathcal G$, their Lie algebras are the same.

Let $G'$ be the monodromy group of $\mathcal G$. The proof of \cite{gkm},11.5.2.4 shows that we have a faithful representation
$G'\hookrightarrow GL(2)$ if $n$ is even and $G'\hookrightarrow
SO(3)\times\mu_n\subset GL(3)$ if $n$ is odd. Let $\mathcal H$ be the corresponding sheaf. As a representation of
the wild inertia group $P_0$ at $0$, the breaks of ${\mathcal G}$ are $0$ and $1/(n+1-m)$, so the breaks of $\mathcal H$ are at most $
1/(n+1-m)$. In particular, the Swan conductor of $\mathcal H$ as a
representation of $P_0$ is $\leq 2/(n+1-m)$ if $n$ is even ($\leq
3/(n+1-m)$ if $n$ is odd). If $n+1-m> 3$ (or $>2$ if $n$ is even), this automatically
implies that $\mathcal H$ is tame at zero as a representation of $\pi_1(\GG_{m,{\bar
k}})$ (since the Swan conductor is an integer) and therefore if factors through the abelian tame fundamental group of $\GG_m$. In particular, the monodromy group would be finite, which contradicts the assumption that $\mathfrak{g}=\mathfrak{sl}_2$. This rules out the possibility $\mathfrak{g}=\mathfrak{sl}_2$ for all cases except $(p,n)=(2,3)$, $(2,5)$ or $(3,2)$.

Therefore the classification theorem forces ${\mathfrak
g}=\mathfrak{sp}_n$ if $n$ is even and ${\mathfrak
g}=\mathfrak{so}_n$ if $n$ is odd as long as $(p,n)\neq(2,3)$, $(2,5)$, $(2,7)$ or $(3,2)$. So in that case $G=Sp(n,\Phi)$ if $n$ is
even, and $G=SO(n,\Phi)$ if $n$ is odd (since the determinant of $\mathcal F$ is geometrically trivial by Proposition \ref{determinant}).

If $(p,n)=(2,3)$, $(2,7)$ or $(3,2)$, $n+1$ is a power of $p$, so $\mathcal F$ is totally wild at $0$ with Swan conductor $1$. By \cite{gkm}, Theorem 8.7.1, applied to the sheaf $\iota^\star{\mathcal F}$ (where $\iota:\GG_m\to\GG_m$ is the inversion map), $\iota^\star{\mathcal F}$ is just a translation of a Kloosterman sheaf on $\GG_m$, so it has the same geometric monodromy group. Using \cite{gkm}, Theorem 11.1, we conclude that $G=Sp(n,\Phi)$ if $(p,n)=(3,2)$, $G=SO(n,\Phi)$ if $(p,n)=(2,3)$ and $G=G_2$ if $(p,n)=(2,7)$.

For the remaining case $p=2$, $n=5$, we have two possibilities, $\mathfrak{g}=\mathfrak{so}_5$ or $\mathfrak{g}=\mathfrak{sl}_2$ in the fourth symmetric power of its standard representation. In the first case, $G$ would be $SO(5)$, since the determinant is trivial. We will rule out this possibility by computing the third moment of ${\mathcal F}$ over ${\mathbb F}_{2^{16}}$. Suppose that $G=SO(5)$, and let $V$ be the stalk of $\mathcal F$ at the generic point of $\AAA^1$, viewed as a representation of $SO(5)$. The altertating square of $\wedge^2 V$ of $V$ is irreducible, and the symmetric square $\mathrm{sym}^2V$ contains the trivial representation and another irreducible factor $W$. So $V\otimes V$ decomposes as $\wedge^2V\oplus{\mathbf 1}\oplus W$. None of these irreducible factors is isomorphic to $V$, so $V\otimes V\otimes V\cong\mathrm{Hom}_G(V\otimes V,V)$ (since $V$ is self-dual) does not contain the trivial representation. Therefore $\HH^2_c(\GG_{m,\bar k},{\mathcal F}^{\otimes 3})$ vanishes, being the dual of $(V\otimes V\otimes V)^G=0$. Since ${\mathcal F}^{\otimes 3}$ does not have punctual sections, its $\HH^0_c$ vanishes too, and then the trace formula gives
$$
\left|\sum_{t\in k^\star}\mathrm{Tr}(F_t|{\mathcal F}_t)^3\right|=|\mathrm{Tr}(F|\HH^1_c(\GG_{m,\bar k},{\mathcal F}^{\otimes 3}))|\leq\dim\HH^1_c(\GG_{m,\bar k},{\mathcal F}^{\otimes 3})q^{6+\frac{1}{2}}
$$
since ${\mathcal F}^{\otimes 3}$ is pure of weight 12. Now ${\mathcal F}^{\otimes 3}$ has rank 125, it is smooth on $\GG_m$, tamely ramified at infinity and all its breaks at $0$ are $\leq 1$ (since the only breaks of $\mathcal F$ at $0$ are $0$ and $1$). Therefore its Swan conductor at $0$ is at most $125$, and then the Euler-Poincar\'e formula gives
$$
\dim H^1_c(\GG_{m,\bar k},{\mathcal F}^{\otimes 3})=-\chi(\GG_{m,\bar k},{\mathcal F}^{\otimes 3})=\mathrm{Sw}_0({\mathcal F}^{\otimes 3})\leq 125
$$
so
$$
\left|\sum_{t\in k^\star}\mathrm{Tr}(F_t|{\mathcal F}_t)^3\right|\leq 125 \cdot q^{6+\frac{1}{2}}.
$$
Now using the explicit formula given in Proposition \ref{gauss}, we find for $k={\mathbb F}_{2^{16}}$ that
$$
\sum_{t\in k^\star}\mathrm{Tr}(F_t|{\mathcal F}_t)^3\simeq 5.48857\cdot 10^{33}>2.5353\cdot 10^{33}\simeq 125\cdot 2^{16(6+\frac{1}{2})}
$$
in contradiction with the inequality above. So $\mathfrak{g}=\mathfrak{sl}_2$ in $\mathrm{sym}^4$ of its standard representation, and therefore $G_0=SL(2)$ in $\mathrm{sym}^4$ of its standard representation. $G_0$ is normal in $G$, being its identity component. For every $g\in G$, conjugation by $g$ gives an automorphism of $G_0$. But every automorphism of $SL(2)$ is inner, so there is an element $g_0\in G_0$ such that $gg_0^{-1}$ is in the centralizer of $G_0$. Now the centralizer of $G_0$ in $GL(5)$ is the set of scalar matrices (a matrix commuting with all matrices of the form
$$
\mathrm{sym}^4\left(
                 \begin{array}{cc}
                   1 & a \\
                   0 & 1 \\
                 \end{array}
               \right)
               \mbox{ and sym}^4\left(
                               \begin{array}{cc}
                                 1 & 0 \\
                                 a & 1 \\
                               \end{array}
                             \right)
                             $$
                             must already be a scalar). But $G\subset SO(5)$, and the only scalar matrix in $SO(5)$ is the identity. Therefore, $g=g_0\in G_0$, and $G=G_0=SL(2)$ in $\mathrm{sym}^4$ of its standard representation. \hfill$\Box$

\section{$L$-functions of symmetric and alternating powers of $\mathcal
F$}

Throughout this section we will assume that $n+1$ is prime to $p$. We will
describe the $L$-function of the smooth sheaf $\Sym^a{\mathcal F}\otimes\wedge^b{\mathcal
F}$ on the set $U=\AAA^1_k-\{(n+1)\zeta:\zeta^{n+1}=1\}$. For simplicity, we will
assume that $k={\mathbb F}_q$, with $(n+1)|(q-1)$, which is always true after a
finite extension of the base field.

\begin{prop} \label{lfunc} The $L$-function of $\mathcal F$ on $U$ is given by
$$
L(U,{\mathcal F},T)=(1-T)P(T)^{n+1}
$$
where $P(T)\in 1+T{\mathbb Z}[T]$ is a polynomial of degree $n-1$.
If $n$ is odd, all reciprocal roots of $P(T)$ have absolute value
$q^{(n-1)/2}$. If $n$ is even, $P(T)=(1\pm q^{(n-2)/2}T)P_1(T)$,
where all reciprocal roots of $P_1(T)$ have absolute value
$q^{(n-1)/2}$.
\end{prop}

{\bf Proof.} Since $\mathcal F$ is smooth, geometrically irreducible and not
geometrically constant on $U$, $L(U,{\mathcal F},T)=\det(1-F\cdot
T|\HH^1_c(U\otimes\bar k,{\mathcal F}))$. If $j:U\to\PP^1$ is the inclusion, the
Euler-Poincar\'e formula gives $\chi(\PP^1_{\bar k},j_\star{\mathcal F})=1+n-(n+1)=0$.
Therefore, $\HH^i(\PP^1_{\bar k},j_\star{\mathcal F})=0$ for all $i$, and we get an
isomorphism
$$
\HH^1_c(U\otimes\bar k,{\mathcal F})\cong(\bigoplus_{\zeta^{n+1}=1}{\mathcal
F}^{I_{(n+1)\zeta}})\oplus{\mathcal F}^{I_\infty}
$$
A similar argument gives
$$
{\mathcal F}^{I_\infty}\cong\HH^1_c(\AAA^1_{\bar k},{\mathcal F}).
$$
By Proposition \ref{lfunction}, we have then
$$
L(U,{\mathcal F},T)=(1-T)\prod_{\zeta^{n+1}=1}\det(1-F\cdot T|{\mathcal
F}^{I_{(n+1)\zeta}})
$$
But the isomorphism ${\mathcal F}\cong [\zeta]^\star{\mathcal F}$ implies
that $P(T)=\det(1-F\cdot T|{\mathcal F}^{I_{(n+1)\zeta}})$ is
independent of $\zeta$. The absolute values of the reciprocal
roots of $P$ are given by Proposition \ref{root}.\hfill$\Box$

\bigskip

We now turn to the study of the $L$-function of the sheaf ${\mathcal
G}_{a,b}:=\Sym^a{\mathcal F}\otimes\wedge^b{\mathcal F}$, which is smooth of rank
${{n+a-1}\choose{a}}{n\choose b}$ and pure of weight $(a+b)(n-1)$ on $U$. Let us
find the bad factor of the $L$-function at infinity first. The local monodromy of ${\mathcal
G}_{a,b}$ it infinity is clearly unipotent, since that of ${\mathcal F}$ is. By Proposition
\ref{lfunc}, the eigenvalues of the geometric Frobenius element at infinity acting on
${\mathcal G}_{a,b}$ are $q^{i_1+\cdots+i_a+j_1+\cdots+j_b}$ for all possible choices of
integers $0\leq i_1\leq i_2\leq\cdots\leq i_a\leq n-1$ and $0\leq
j_1<j_2<\cdots<j_b\leq n-1$. Let $N_{n,a,b,k}$ be the number of such possible
choices with $i_1+\cdots+i_a+j_1+\cdots+j_b=k$, that is,
$$
N_{n,a,b,k}=\#\{(i_1,\ldots,i_a,j_1,\ldots,j_b):0\leq i_1\leq
i_2\leq\cdots\leq i_a\leq n-1, $$ $$ 0\leq j_1<j_2<\cdots<j_b\leq
n-1,i_1+\cdots+i_a+j_1+\cdots+j_b=k\}
$$
It is clear that $N_{n,a,b,k}=N_{n,a,b,(a+b)(n-1)-k}$ (just change
$i_l\mapsto n-1-i_{a+1-l}$ and $j_l\mapsto n-1-j_{b+1-l}$) and
$N_{n,a,b,k}=0$ for $k<b(b-1)/2$ and $k>(a+b)(n-1)-b(b-1)/2$.

\begin{prop} The dimension of the invariant subspace ${\mathcal
G}_{a,b}^{I_\infty}$ is $N_{n,a,b,c}$ where
$c=\lfloor\frac{(a+b)(n-1)}{2}\rfloor$ and $N_{n,a,b,c}$ is the
coefficient of $x^cz^b$ in the expansion of the power series
$$
\frac{(1-x^n)\cdots(1-x^{a+n-1})}{(1-x)\cdots(1-x^a)}(1+z)(1+xz)\cdots(1+x^{n-1}z).
$$
If $(a+b)(n+1)$ is even, all Jordan blocks for the action of
$I_\infty$ on ${\mathcal G}_{a,b}$ have odd size, and the number of
blocks of size $2k+1$ is $N_{n,a,b,c-k}-N_{n,a,b,c-k-1}$ for all
$k\geq 0$. If $(a+b)(n+1)$ is odd, all Jordan blocks for the
action of $I_\infty$ on ${\mathcal G}_{a,b}$ have even size, and the
number of blocks of size $2k+2$ is $N_{n,a,b,c-k}-N_{n,a,b,c-k-1}$
for all $k\geq 0$.
\end{prop}

{\bf Proof.} This is just a translation of \cite{weilii}, 1.8.4
and \cite{gkm}, 7.0.7 to this particular situation, considering
that ${\mathcal G}_{a,b}$ is pure of weight $(a+b)(n-1)$ and all
Frobenius eigenvalues of ${\mathcal G}_{a,b}$ at infinity are integral
powers of $q$ (that is, they have even weight). In fact, the
multiplicity $N_{n,a,b,0}$ of the minimun Froebnius eigenvalue
$q^0$ is equal to the number of Jordan blocks with length
$(a+b)(n-1)+1$. Removing these blocks, then the multiplicity
$N_{n,a,b,1}-N_{n,a,b,0}$ of the minimun remaining Frobenius
eigenvalue $q$ is equal to the number of blocks with length
$(a+b)(n-1)-1$. By induction, for $0<k\leq c$, one deduces that
$N_{n,a,b,k}-N_{n,a,b,k-1}$ is equal to the number of blocks with
length $(a+b)(n-1)-2k+1$ and with minimun Frobenius eigenvalue
$q^k$. The dimension of the invariant subspace ${\mathcal
G}_{a,b}^{I_\infty}$ is simply the total number of Jordan blocks:
$$\sum_{k=0}^c (N_{n,a,b,k}-N_{n,a,b,k-1})=   N_{n,a,b,c}.$$
\hfill$\Box$

\begin{cor} The local $L$-function of $j_\star{\mathcal G}_{a,b}$ at infinity
has degree $N_{n,a,b,c}$ and is given by
$$
\det(1-F_\infty\cdot T|{\mathcal
G}_{a,b}^{I_\infty})=\prod_{k=0}^c(1-q^kT)^{\alpha_{a,b}(k)}
$$
where $\alpha_{a,b}(k)=N_{n,a,b,k}-N_{n,a,b,k-1}$.\end{cor}

\bigskip

We can construct a generating function for $\alpha(k)$ in the following way. Let
$C_{n,a,k}=\#\{(i_1,\ldots,i_a):0\leq i_1\leq i_2\leq\cdots\leq i_a\leq
n-1,i_1+\cdots+i_a=k\}=\#\{(h_0,\ldots,h_{n-1}):0\leq
h_i,h_0+\cdots+h_{n-1}=a,h_1+2h_2+\cdots+(n-1)h_{n-1}=k\}$ (to check that both
numbers agree, just let $h_j$ be the number of $l=1,\ldots,a$ such that $i_l=j$). By
\cite{fw2}, Theorem 3.1, we have
$$
\sum_{k\geq
0}(C_{n,a,k}-C_{n,a,k-1})x^k=\{\frac{(1-x^n)\cdots(1-x^{n+a-1})}{(1-x^2)\cdots(1-x^a)}\},
$$
where the quantity in the bracket is understood to be $1-x^n$ if
$a=1$, and $1-x$ if $a=0$. Let
$$B_{n,b,j}=\#\{(j_1,\ldots,j_b):0\leq j_1<\cdots<j_b\leq
n-1,j_1+\cdots+j_b=j\}.$$It is the coefficient of $x^jz^b$ in the
expansion of $(1+z)(1+xz)\cdots(1+x^{n-1}z)$. Then
$$
N_{n,a,b,k}=\sum_{j=0}^k C_{n,a,k-j}B_{n,b,j},
$$
and thus
$$
\alpha_{a,b}(k)=N_{n,a,b,k}-N_{n,a,b,k-1}=\sum_{j=0}^{k-1}(C_{n,a,k-j}-C_{n,a,k-j-1})B_{n,b,j}+B_{n,b,k}.
$$
Therefore $\alpha_{a,b}(k)$ is the coefficient of $x^kz^b$ in the
expansion of
$$
\{\frac{(1-x^n)\cdots(1-x^{a+n-1})}{(1-x^2)\cdots(1-x^a)}\}(1+z)(1+xz)\cdots(1+x^{n-1}z).
$$
In particular, the number $N_{n,a,b,c}$ is the coefficient of
$x^cz^b$ in the expansion of the power series
$$
\frac{(1-x^n)\cdots(1-x^{a+n-1})}{(1-x)\cdots(1-x^a)}(1+z)(1+xz)\cdots(1+x^{n-1}z).
$$

We now look for the bad factors of the $L$-function at the finite
singular points $t=(n+1)\zeta$ with $\zeta^{n+1}=1$. Suppose that
$n$ is even. Then the local monodromy at $t$ is unipotent, with a
Jordan block of size $2$ and all other blocks of size $1$. The
Frobenius eigenvalues on ${\mathcal F}^{I_t}$ are $\epsilon
q^{(n-2)/2}$, with $\epsilon=1$ or $-1$, and $(n-2)/2$ pairs of
conjugate complex numbers
$\alpha_1,\ldots,\alpha_{(n-1)/2},\bar\alpha_1,\ldots,\bar\alpha_{(n-1)/2}$
of absolute value $q^{(n-1)/2}$. That is, as a representation of
$I_t$, ${\mathcal F}\cong U_2\oplus{\mathbf 1}^{n-2}$, where $U_m$
denotes the unique (up to isomorphism) non-trivial unipotent tame
representation of $I_t$ of dimension $m$ with a single Jordan
block. Therefore, we get isomorphisms
$$
\Sym^a{\mathcal F}\cong \bigoplus_{i=0}^a \Sym^i
U_2\otimes\Sym^{a-i}{\mathbf 1}^{n-2}=\bigoplus_{i=0}^a
U_{i+1}^{{n-3+a-i}\choose{n-3}}
$$
$$
\wedge^b{\mathcal F}\cong \wedge^b{{\mathbf
1}^{n-2}}\oplus(U_2\otimes\wedge^{b-1}{{\mathbf 1}^{n-2}})\oplus
\wedge^{b-2}{{\mathbf 1}^{n-2}}\cong {\mathbf 1}^{{{n-2}\choose
b-2}+{{n-2}\choose b}}\oplus U_2^{{n-2}\choose {b-1}}.
$$

\begin{lem} Let $V$ and $W$ be vector spaces of dimensions $n\geq 2$ and
$2$ respectively over an algebraically closed field $k$ of
characteristic $0$, and let $T:V\to V$ and $U:W\to W$ be unipotent
endomorphisms with a single Jordan block. Then $T\otimes U:V\otimes
W\to V\otimes W$ is unipotent with two Jordan blocks of sizes $n+1$
and $n-1$.
\end{lem}

{\bf Proof.} Let $\{\xx,\yy\}$ be a basis for $W$ such that
$U(\xx)=\xx$ and $U(\yy)=\xx+\yy$. We claim that the invariant
subspace of $T\otimes U$ is the subspace of elements that can be
written as $\vv\otimes\xx+(\vv-T(\vv))\otimes\yy$ for
$\vv\in{\mathrm{Ker}}((T-I_V)^2)$, which has dimension $2$ by
hypothesis:
$$
(T\otimes
U)(\vv\otimes\xx+(\vv-T(\vv))\otimes\yy)=T(\vv)\otimes\xx+T(\vv-T(\vv))\otimes(\xx+\yy)$$
$$
=T(\vv)\otimes\xx+(\vv-T(\vv))\otimes(\xx+\yy)=\vv\otimes\xx+(\vv-T(\vv))\otimes\yy.
$$
Conversely, if $(T\otimes
U)(\vv\otimes\xx+\ww\otimes\yy)=\vv\otimes\xx+\ww\otimes\yy$, we get
$T(\ww)=\ww$ and $T(\vv)+T(\ww)=\vv$, so $\ww=T(\ww)=\vv-T(\vv)$ and
$(T-I_V)^2(\vv)=0$. This shows that $T\otimes U$ has precisely two
Jordan blocks. From
$$
T\otimes U-I\otimes I=(T-I)\otimes(U-I)+I\otimes(U-I)+(T-I)\otimes I
$$
we get that $(T\otimes U-I\otimes I)^{n+1}$ is a sum of terms
$(T-I)^\alpha\otimes(U-I)^\beta$ with $\alpha+\beta\geq n+1$ and
therefore equal to $0$, since $(T-I)^n=(U-I)^2=0$. So the Jordan
blocks of $T\otimes U$ have size $\leq n+1$. Finally, if $\vv\in V$
is a vector such that $\ww:=(T-I)^{n-1}(\vv)\neq 0$ and $\xx,\yy\in
W$ are as above, the same expression shows that $(T\otimes
U-I\otimes
I)^n(\vv\otimes\yy)=(n-1)(T-I)^{n-1}(\vv)\otimes(U-I)(\yy)=(n-1)\ww\otimes\xx\neq
0$, so $\vv\otimes\yy$ generates a Jordan block of size $n+1$, and
the other block must have size $2n-(n+1)=n-1$.\hfill$\Box$

\begin{cor} Suppose that $n$ is even. As a representation of $I_t$, ${\mathcal G}_{a,b}=\Sym^a{\mathcal
F}\otimes\wedge^b{\mathcal F}$ is isomorphic to
$$
\bigoplus_{i=0}^a
U_{i+1}^{{{n-3+a-i}\choose{n-3}}\left[{{n-2}\choose{b-2}}+{{n-2}\choose
b}\right]}\oplus
U_i^{{{n-3+a-i}\choose{n-3}}{{n-2}\choose{b-1}}}\oplus
U_{i+2}^{{{n-3+a-i}\choose{n-3}}{{n-2}\choose{b-1}}}=\bigoplus_{i=1}^{a+2} U_i^{d(i)}
$$
where
$d(i)=\left[{{n-3+a-i}\choose{n-3}}+{{n-1+a-i}\choose{n-3}}\right]{{n-2}\choose{b-1}}+{{n-2+a-i}\choose{n-3}}\left[{{n-2}\choose{b-2}}+{{n-2}\choose
b}\right]$.\end{cor}

\begin{cor} \label{localeven} Suppose that $n$ is even. The local $L$-function of $j_\star{\mathcal G}_{a,b}$ at $t$,
$\det(1-F_t\cdot T|{\mathcal G}_{a,b}^{I_t})$ has degree \bigskip

$D_{n,a,b}:=
\sum_{i=1}^{a+2}d(i)=\left[{{n-3+a}\choose{n-2}}+{{n-1+a}\choose{n-2}}\right]{{n-2}\choose{b-1}}+{{n-2+a}\choose{n-2}}\left[{{n-2}\choose{b-2}}+{{n-2}\choose
b}\right]. $

\bigskip
For every $i=1,\ldots,a+2$, it has $d(i)$
roots which are pure of weight $(a+b)(n-1)-(i-1)$.
\end{cor}

\bigskip

For $n$ odd, the situation is much simpler. In that case, as a
representation of $I_t$, ${\mathcal F}\cong \chi_2\oplus{\mathbf
1}^{n-1}$, where $\chi_2:I_t\to\QQ^\star$ is the unique character of
order $2$. Therefore, we get isomorphisms
$$
\Sym^a{\mathcal F}\cong \bigoplus_{i=0}^a \Sym^i
\chi_2\otimes\Sym^{a-i}{\mathbf 1}^{n-1}\cong\bigoplus_{i=0 \atop
{i\mathrm{even}}}^a {\mathbf
1}^{{n-2+a-i}\choose{n-2}}\oplus\bigoplus_{i=0
\atop{i\mathrm{odd}}}^a \chi_2^{{n-2+a-i}\choose{n-2}}
$$
$$
\wedge^b{\mathcal F}\cong \wedge^b{{\mathbf
1}^{n-1}}\oplus(\chi_2\otimes\wedge^{b-1}{{\mathbf 1}^{n-1}})\cong
{\mathbf 1}^{{n-1}\choose b}\oplus \chi_2^{{n-1}\choose {b-1}}.
$$
$$
\Sym^a{\mathcal F}\otimes\wedge^b{\mathcal F}\cong {\mathbf
1}^{\alpha{{n-1}\choose
b}+\beta{{n-1}\choose{b-1}}}\oplus\chi_2^{\alpha{{n-1}\choose
{b-1}}+\beta{{n-1}\choose{b}}}
$$
where
$$
\alpha=\sum_{i=0 \atop {i\mathrm{even}}}^a {{n-2+a-i}\choose{n-2}},
\beta=\sum_{i=0 \atop {i\mathrm{odd}}}^a {{n-2+a-i}\choose{n-2}}
$$

\begin{cor} \label{localodd} Suppose that $n$ is odd. The local $L$-function of $j_\star{\mathcal G}_{a,b}$ at $t$,
$\det(1-F_t\cdot T|{\mathcal G}_{a,b}^{I_t})$ has degree
$$D_{n,a,b}:=
{{n-1}\choose b}\sum_{i=0 \atop {i\mathrm{even}}}^a
{{n-2+a-i}\choose{n-2}}+{{n-1}\choose{b-1}}\sum_{i=0 \atop
{i\mathrm{odd}}}^a {{n-2+a-i}\choose{n-2}}.
$$
All its roots are pure of weight $(a+b)(n-1)$.
\end{cor}

Consider the sheaf $j_\star{\mathcal G}_{a,b}$ on $\PP^1$. The
Tate-twisted sheaf ${\mathcal G}_{a,b}((n-1)(a+b)/2)$ is self-dual, so
Poincar\'e duality gives a perfect pairing of $\Gal(\bar
k/k)$-modules
$$
H^i(\PP^1_{\bar k},j_\star{\mathcal G}_{a,b})\times H^{2-i}(\PP^1_{\bar
k},j_\star{\mathcal G}_{a,b})\to \QQ((a+b)(1-n)-1)
$$
for $i=0,1,2$. Since ${\mathcal G}_{a,b}$ is smooth on $U$, the zeroth cohomology
group $H^0(\PP^1_{\bar k},j_\star{\mathcal G}_{a,b})$ corresponds to the maximal
geometrically constant subsheaf of ${\mathcal G}_{a,b}$. Since ${\mathcal G}_{a,b}$ is pure
of weight $(n-1)(a+b)$ and all Frobenius eigenvalues of $j_\star{\mathcal G}_{a,b}$ at
infinity are integral powers of $q$, such a subsheaf must be a direct sum of copies
of $\QQ((1-n)(a+b)/2)$. Incidentally, this shows that $H^0(\PP^1_{\bar k},j_\star{\mathcal
G}_{a,b})=0$ if $(n-1)(a+b)$ is odd. Therefore, we have:

\begin{prop} The $L$-function of $j_\star{\mathcal G}_{a,b}$ on $\PP^1$
has the form
$$
L(\PP^1,j_\star{\mathcal
G}_{a,b})=\frac{P_{a,b}(T)}{(1-q^{(a+b)(n-1)/2}T)^{\delta_{a,b}}(1-q^{(a+b)(n-1)/2+1}T)^{\delta_{a,b}}}
$$
where $\delta_{a,b}=\dim H^0(\PP^1_{\bar k},j_\star{\mathcal G}_{a,b})$,
and $P_{a,b}(T)$ is a polynomial that satisfies the functional
equation
$$
P_{a,b}(T)=\pm T^r q^{((a+b)(n-1)+1)r/2}
P_{a,b}(1/q^{(a+b)(n-1)+1}T)
$$
where $r=\deg(P_{a,b})$.\end{prop}

{\bf Proof.} We have just seen that $H^0(\PP^1_{\bar
k},j_\star{\mathcal G}_{a,b})=\QQ((1-n)(a+b)/2)^{\delta_{a,b}}$, and
Poincar\'e duality implies that $H^2(\PP^1_{\bar k},j_\star{\mathcal
G}_{a,b})=\QQ((1-n)(a+b)/2-1)^{\delta_{a,b}}$. This gives the
denominator.

 The numerator is
$P_{a,b}(T)=(1-\alpha_1T)\cdots(1-\alpha_rT)$, where
$\alpha_1,\ldots,\alpha_r$ are the Frobenius eigenvalues of
$H^1(\PP^1_{\bar k},j_\star{\mathcal G}_{a,b})$. By Poincar\'e
duality, these eigenvalues are permuted by $\alpha\mapsto
q^{(a+b)(n-1)+1}/\alpha$. In particular,
$(\prod\alpha_i)^2=q^{((a+b)(n-1)+1)r}$. We have
$$
P_{a,b}(1/q^{(a+b)(n-1)+1}T)=(1-\frac{1}{\alpha_1T})\cdots(1-\frac{1}{\alpha_rT})$$
$$
=\frac{1}{\alpha_1\cdots\alpha_r
T^r}(\alpha_1T-1)\cdots(\alpha_rT-1)=\frac{(-1)^r}{\pm T^r
q^{((a+b)(n-1)+1)r/2}} P_{a,b}(T)
$$
and the functional equation follows.\hfill$\Box$

\bigskip

To find the dimension of $H^0(\PP^1_{\bar k},j_\star{\mathcal G}_{a,b})$ we will use the
knowledge of the global monodromy of $\mathcal F$, as in \cite{fs}. Let $V$ be the
geometric generic fibre of $\mathcal F$, regarded as a representation of
$\pi_1(U\otimes\bar k)$. We know that the Zariski closure $G$ of the image of
$\pi_1(U\otimes\bar k)$ in $\mathrm{GL}(V)$ is ${Sp}(n)$ if $n$ is even and
${O}(n)$ if $n$ is odd. The dimension we are looking for is the dimension of the
invariant subspace
$\dim(\Sym^a(V)\otimes\wedge^b(V))^G=\dim\mathrm{Hom}_G(\Sym^a(V),\wedge^b(V))$
(since $V$ is self-dual as a representation of $G$).

Suppose $n=2m$ is even. The representations of $G=Sp(n)$ are in one to one
correspondence with the representations of the Lie algebra
$\mathfrak{g}=\mathfrak{sp}_n$. If $L_1,\ldots,L_m$ are generators of the weight
lattice for $\mathfrak g$, then $\mathrm{Sym}^d V$ is the irreducible representation
with maximal weight $dL_1$, and the kernel of the natural contraction map
$\wedge^d V\to \wedge^{d-2} V$ is the irreducible representation of maximal weight
$L_1+\ldots+L_d$ for $1\leq d\leq m$. (\cite{fh}, ch.17) Therefore we have
$$
\wedge^b V\cong W(L_1+\ldots+L_b)\oplus
W(L_1+\ldots+L_{b-2})\oplus\ldots\oplus V
$$
if $b\leq m$ is odd and
$$
\wedge^b V\cong W(L_1+\ldots+L_b)\oplus
W(L_1+\ldots+L_{b-2})\oplus\ldots\oplus {\mathbf 1}
$$
if $b\leq m$ is even and $\wedge^b V\cong \wedge^{n-b} V$ for $m\leq
b\leq n$. So $\mathrm{Sym}^a V\otimes \wedge^b V$ contains exactly
one copy of the trivial representation if $a=0$ and $b\leq n$ is
even or if $a=1$ and $b\leq n$ is odd, and does not contain the
trivial representation otherwise.

Suppose $n=2m+1$ is odd. The representations of $SO(n)$, the
connected component of $G$ containing the identity, are in
one-to-one correspondence with the representations of the Lie
algebra $\mathfrak{g}=\mathfrak{so}_n$ contained in the tensor
algebra of the standard representation. Each of them gives rise to
two different representations of $O(n)$ (given one of them, the
other one is obtained by tensoring with the determinant). If
$L_1,\ldots,L_m$ are generators of the weight lattice for $\mathfrak
g$, then $\wedge^d V$ is the irreducible representation with maximal
weight $L_1+\ldots+L_d$ for $d\leq m$, $\wedge^d
V\cong\wedge^{n-d}V$ for $m+1\leq d\leq n$, and the kernel of the
natural contraction map $\mathrm{Sym}^d V\to \mathrm{Sym}^{d-2} V$
is the irreducible representation of maximal weight $dL_1$
(cf. \cite{fh}, ch.19). Therefore we have
$$\mathrm{Sym}^a V\cong W(aL_1)\oplus
W((a-2)L_1)\oplus\ldots\oplus V
$$
if $a$ is odd and
$$
\mathrm{Sym}^a V\cong W(aL_1)\oplus W((a-2)L_1)\oplus\ldots\oplus
{\mathbf 1}
$$
if $a$ is even. So $\mathrm{Sym}^a V\otimes \wedge^b V$ (as a
representation of $\mathfrak g$) contains exactly one copy of the
trivial representation if $a$ is even and $b=0$ or $n$, or if $a$ is
odd and $b=1$ or $n-1$, and does not contain the trivial
representation otherwise.

For $G$ itself, since the determinant becomes trivial only in even tensor powers of
the standard representation, we get that $\mathrm{Sym}^a V\otimes \wedge^b V$
contains exactly one copy of the trivial representation and no copies of the
determinant representation if $a$ is even and $b=0$, or if $a$ is odd and $b=1$, it
contains exactly one copy of the determinant representation and no copies of the
trivial representation if $a$ is even and $b=n$ or if $a$ is odd and $b=n-1$, and it
does not contain the trivial or the determinant representations otherwise. Therefore
we get:

\begin{prop}\label{dimension} The dimension $\delta_{a,b}=\dim H^0(\PP^1_{\bar k},j_\star{\mathcal
G}_{a,b})$ is
$$
\begin{array}{ll}
\mbox{if $n$ is even } & \left\{\begin{array}{ll} 1 & \mbox{if $a=0$
and
$b\leq n$ is even or $a=1$ and $b\leq n$ is odd} \\
0 & \mbox{otherwise} \end{array}\right.
\\
\mbox{if $n$ is odd } & \left\{\begin{array}{ll} 1 & \mbox{if $a$ is
even and
$b=0$ or $a$ is odd and $b=1$} \\
0 & \mbox{otherwise} \end{array}\right.
\end{array}
$$
\end{prop}

Putting everything together, we get the following expression for the $L$-function of
${\mathcal G}_{a,b}$:

\begin{thm} The $L$-function of ${\mathcal G}_{a,b}$ on $U$ has
total degree $n{{n+a-1}\choose{a}}{n\choose b}$ and is given
by
$$
L(U,{\mathcal
G}_{a,b})=\frac{P_{a,b}(T)Q_{a,b}(T)^{n+1}\prod_{k=0}^{[(a+b)(n-1)/2]}(1-q^kT)^{\alpha_{a,b}(k)}}
{(1-q^{(a+b)(n-1)/2}T)^{\delta_{a,b}}(1-q^{(a+b)(n-1)/2+1}T)^{\delta_{a,b}}}
$$
where $\delta_{a,b}=0$ or $1$ is given by Proposition
\ref{dimension}, $\alpha_{a,b}(k)=N_{n,a,b,k}-N_{n,a,b,k-1}$,
$Q_{a,b}(T)$ is a polynomial whose degree $D_{n,a,b}$ and the
weights of its roots are given in Corollaries \ref{localeven} and
\ref{localodd} and $P_{a,b}(T)$ is a polynomial in $1+T\Z[T]$ of
degree
$$n{n+a-1\choose a}{n\choose
b}+2\delta_{a,b}-N_{n,a,b,c}-(n+1)D_{n,a,b},$$ where
$c=[\frac{(a+b)(n-1}{2}]$. Furthermore, $P_{a,b}(T)$ is pure of
weight $(a+b)(n-1)+1$ and it satisfies the functional equation
$$
P_{a,b}(T)=\pm T^r q^{((a+b)(n-1)+1)r/2}
P_{a,b}(1/q^{(a+b)(n-1)+1}T).
$$
\end{thm}

{\bf Proof.} The total degree of the $L$-function is the negative Euler
characteristic $-\chi(U,{\mathcal G}_{a,b})$. Since ${\mathcal G}_{a,b}$ is
everywhere tamely ramified, this Euler characteristic is
$\chi(U)\mathrm{rank}({\mathcal G}_{a,b})=-n{{n+a-1}\choose{a}}{n\choose
b}$. The stated formula is just the decomposition
$$
L(U,{\mathcal G}_{a,b})=L(\PP^1,j_\star{\mathcal
G}_{a,b})\det(1-F_\infty\cdot T|{\mathcal
G}_{a,b}^{I_\infty})\prod_{t\in(n+1)\mu_{n+1}}\det(1-F_t\cdot
T|{\mathcal G}_{a,b}^{I_t}),
$$
and the shape of each of the factors has already been
determined.\hfill$\Box$

\bigskip

The fact that
$$L({\mathbb A}^1, {\mathcal G}_{a,b}) = \frac{L(U,{\mathcal
G}_{a,b})}{Q_{a,b}(T)^{n+1}}$$ together with the above theorem
immediately implies Theorem 1.1.

\begin{cor} The $L$-function of $[{\mathcal F}]^d$ on ${\mathbb A}^1$ is given by
$$
L({\mathbb A}^1,[{\mathcal
F}]^d)=P_d(T)(1-q^{\frac{d(n-1)}{2}}T)(1-q^{\frac{d(n-1)}{2}+1}T)
\prod_{k=0}^{[\frac{n-2}{2}]}\frac{1-q^{dk}T}{1-q^{dk+1}T}
$$
if $n$ and $d$ are even,
$$
L({\mathbb A}^1,[{\mathcal
F}]^d)=P_d(T)\prod_{k=0}^{[\frac{n-2}{2}]}\frac{1-q^{dk}T}{1-q^{dk+1}T}
$$
if $n$ is even and $d$ is odd,
$$
L({\mathbb A}^1,[{\mathcal
F}]^d)={P_d(T)}(1-q^{\frac{d(n-1)}{2}+1}T)^{-1}\prod_{k=0}^{[\frac{n-2}{2}]}\frac{1-q^{dk}T}{1-q^{dk+1}T}
$$
if $n$ is odd and $d$ is even and
$$
L({\mathbb A}^1,[{\mathcal
F}]^d)=P_d(T)(1-q^{\frac{d(n-1)}{2}}T)\prod_{k=0}^{[\frac{n-2}{2}]}\frac{1-q^{dk}T}{1-q^{dk+1}T}
$$
if $n$ and $d$ are odd, where
$$
P_d(T)=\prod_{b=0}^n P_{d-b,b}(T)^{(-1)^{b-1}(b-1)}.
$$
Alternatively, the above four expressions can be unified into the
following single expression
$$L({\mathbb A}^1,[{\mathcal
F}]^d)=P_d(T)\frac{
(1-q^{\frac{d(n-1)}{2}}T)^{\frac{1+(-1)^{d+n}}{2}}}
{(1-q^{\frac{d(n-1)}{2}+1}T)^{\frac{-(-1)^n-(-1)^{n+d}}{2}}}
\prod_{k=0}^{[\frac{n-2}{2}]}\frac{1-q^{dk}T}{1-q^{dk+1}T}.$$
\end{cor}

{\bf Proof.} From the $L$-function decomposition
$$
L({\mathbb A}^1, [{\mathcal F}]^d) = \prod_{b=0}^n L({\mathbb A}^1, {\mathcal G}_{d-b, b})^{(-1)^{b-1}(b-1)}$$
and Theorem 1.1 we get
$$
L({\mathbb A}^1, [{\mathcal F}]^d) = P_d(T)\frac{\prod_{k=0}^{[d(n-1)/2]}(1-q^kT)^{\sum_{b=0}^n (-1)^{b-1}b\alpha_{d-b,b}(k)}}{(1-q^{d(n-1)/2}T)^{\delta_d}(1-q^{d(n-1)/2+1}T)^{\delta_d}}
$$
where $\delta_d=\sum_{b=0}^n (-1)^{b-1}(b-1)\delta_{d-b,b}$. Using
Proposition \ref{dimension}, we find $\delta_d=-d+(d-1)=-1$ if $n$
and $d$ are even, $\delta_d=1$ if $n$ is odd and $d$ is even and
$\delta_d=0$ if $d$ is odd. It remains to compute the numerator of
the previous expression, which is just the local $L$-function at
infinity of the virtual sheaf $[{\mathcal F}]^d$. Write $[{\mathcal
F}]^{d}=[{\mathcal H}_{+}]-[{\mathcal H}_{-}]$, where ${\mathcal
H}_{+}=\oplus_{b=0,b\;\mathrm{ odd}}^n{\mathcal G}_{d-b,b}^b$ and
${\mathcal H}_{-}=\oplus_{b=0,b\;\mathrm{ even}}^n{\mathcal G}_{d-b,b}^b$
are ``real" sheaves. We know that ${\mathcal H}_{+}$ and ${\mathcal
H}_{-}$ are pure of weight ${n-1}$, the inertia group $I_{\infty}$
acts unipotently on them and all their Frobenius eigenvalues at
infinity are integral powers of $q$. If $\mu(k)$ (resp. $\nu(k)$)
is the number of Frobenius eigenvalues of ${\mathcal H}_{+}$ (resp. of
${\mathcal H}_{-}$) at infinity which are equal to $q^k$, the local
factor at infinity of the $L$-function of $[{\mathcal F}]^d$ is given
by $\prod_{k= 0}^{[d(n-1)/2]}
(1-q^kT)^{(\mu(k)-\nu(k))-(\mu(k-1)-\nu(k-1))}$, again by
\cite{weilii}, 1.8.4 and \cite{gkm}, 7.0.7.

On the other hand, for every $r\geq 1$ the trace of the action of the $dr$-th power of the local Frobenius at infinity on $[{\mathcal F}]^d$ is
$$
\Trace(F_{\infty}^{dr}|{\mathcal F})=1+q^{dr}+\cdots+q^{dr(n-1)}.
$$
But
$$
\Trace(F_{\infty}^{dr}|{\mathcal F})=\Trace(F_{\infty}^r|[{\mathcal F}]^d)=\sum_{k\geq 0}(\mu(k)-\nu(k))q^{kr}.
$$
Since this holds for every $r\geq 1$, we conclude that $\mu(k)-\nu(k)=1$ if $k=0,d,\ldots,(n-1)d$ and $0$ otherwise. Therefore, the local factor at infinity of the $L$-function of $[{\mathcal F}]^d$ is
$$
\prod_{k= 0}^{[d(n-1)/2]} (1-q^kT)^{(\mu(k)-\nu(k))-(\mu(k-1)-\nu(k-1))}$$
$$
=(1-q^{d(n-1)/2}T)\prod_{k=0}^{\frac{n-3}{2}}\frac{1-q^{dk}T}{1-q^{dk+1}T}
$$
if $n$ is odd and
$$
\prod_{k= 0}^{[d(n-1)/2]} (1-q^kT)^{(\mu(k)-\nu(k))-(\mu(k-1)-\nu(k-1))}$$
$$
=\prod_{k=0}^{\frac{n}{2}-1}\frac{1-q^{dk}T}{1-q^{dk+1}T}
$$
if $n$ is even. This combined with the explicit description of $\delta_d$ proves the result.\hfill$\Box$

\bigskip

We can now finish the proof of Theorem 1.1. By Theorem 2.1, we
deduce
$$
L({\mathbb A}^1, [{\mathcal H}^{n-1}(K)]^d) =L({\mathbb A}^1, [{\mathcal
F}]^d)L({\mathbb A}^1, \QQ^n)=L({\mathbb A}^1, [{\mathcal
F}]^d)(1-qT)^{-n},$$ and for $n\leq j \leq 2(n-1)$,
$$L({\mathbb A}^1, [{\mathcal H}^{j}(K)]^d)=
L({\mathbb A}^1, \QQ(d(n-1-j))^{{n\choose
j-n+2}}=(1-q^{d(j-(n-1))+1}T)^{-{n \choose j-n+2}}.$$ Also, by
Theorem 2.1 and the Grothendieck trace formula,
$$Z_d({\mathbb A}^1, X_{\lambda}) =\prod_{j=n-1}^{2(n-1)} L({\mathbb A}^1, [{\mathcal
H}^{j}(K)]^d)^{(-1)^j}.$$ Substituting the above calculation, we
obtain
$$Z_d({\mathbb A}^1, X_{\lambda})^{(-1)^{n-1}}=L({\mathbb A}^1, [{\mathcal
F}]^d)\prod_{i=0}^{n-1}(1-q^{di+1}T)^{(-1)^{i+1}{n\choose i+1}}.$$
This together with Corollary 3.11 gives Theorem 1.1. The proof is
complete.

\section{Zeta function in terms of Gauss
sums}

In this section,  we give an elementary formula for the number
$N_q(\lambda)$ of ${\mathbb F}_q$-rational points  in the fibre
$X_\lambda$ in terms of Gauss sums for every $\lambda\in{\mathbb
F}_q$. This type of elementary formulas for a general equation can be found in
Koblitz \cite{ko}. We derive a more explicit formula in the special case of $X_{\lambda}$
and in particular deduce an explicit formula for the zeta
function of $X_0$. This allows us to determine the rank of the sheaf ${\mathcal F}$
when $p$ divides $n+1$ and the local factor at $0$ of the sheaf ${\mathcal F}$.

Let $\omega:{\mathbb F}_q^\star\to {\mathbb C}^*$ be a primitive
character of order $q-1$. For every $k\in{\mathbb Z}$, define the
Gauss sum $G_q(k)$ by
$$
G_q(k)=-\sum_{a\in{\mathbb F}_q^\star}
\omega(a)^{-k}\zeta_p^{\Tr_{{\mathbb F}_q/{\mathbb F}_p}(a)}
$$
where $\zeta_p=\exp(2\pi i/p)$. It is clear that $G_q(k)=1$ if
$(q-1)|k$, and $|G_k(q)|=\sqrt{q}$ otherwise. We have the inversion
formula
\begin{equation} \label{inv}
\zeta_p^{\Tr(a)}=\sum_{k=0}^{q-2}\frac{G_q(k)}{1-q}\omega(a)^k
\end{equation}
for every $a\in{\mathbb F}_q^\star$. We find that
$$
N_q(\lambda)=\frac{1}{q}\sum_{x_0\in\Fq}\sum_{x_1,\ldots,x_n\in\Fq^\star}\zeta_p^{\Tr(x_0x_1+\cdots+x_0x_n+\frac{x_0}{x_1\cdots
x_n}-x_0\lambda)}=\frac{(q-1)^n}{q}+\frac{1}{q}S_q(\lambda)
$$
where
$$
S_q(\lambda)=\sum_{x_0,x_1,\ldots,x_n\in\Fq^\star}\zeta_p^{\Tr(x_0x_1+\cdots+x_0x_n+\frac{x_0}{x_1\cdots
x_n}-x_0\lambda)}.
$$
Using equation \ref{inv}, we deduce for $\lambda\neq 0$:
$$
S_q(\lambda)=\sum_{x_0,x_1,\ldots,x_n\in\Fq^\star}\zeta_p^{\Tr(x_0x_1)}\cdots\zeta_p^{\Tr(x_0x_n)}\zeta_p^{\Tr(x_0/x_1\cdots
x_n)}\zeta_p^{\Tr(-x_0\lambda)}$$
$$
=\sum_{k_1,\ldots,k_{n+2}=0}^{q-2}\frac{G_q(k_1)\cdots
G_q(k_{n+2})}{(1-q)^{n+2}}\sum_{y_i^{q-1}=1}(y_0y_1)^{k_1}\cdots
(y_0y_n)^{k_n}(\frac{y_0}{y_1\cdots
y_n})^{k_{n+1}}(y_0\omega(-\lambda))^{k_{n+2}}
$$
$$
=\sum_{k_1,\ldots,k_{n+2}=0}^{q-2}\frac{G_q(k_1)\cdots
G_q(k_{n+2})}{(1-q)^{n+2}}\omega(-\lambda)^{k_{n+2}}\sum_{y_i^{q-1}=1}y_0^{k_1+\cdots+k_{n+2}}y_1^{k_1-k_{n+1}}\cdots
y_n^{k_n-k_{n+1}}
$$
$$
=(-1)^n\sum_{{a,b=0} \atop (n+1)a+b\equiv 0
(q-1)}^{q-2}\frac{G_q(a)^{n+1}G_q(b)}{q-1}\omega(-\lambda)^b$$
$$
={(-1)^n}\left( \frac{1}{q-1}+\sum_{(n+1)a+b\equiv 0 (q-1) \atop
(a,b)\neq(0,0)} \frac{G_q(a)^{n+1}
G_q(b)}{q-1}\omega(-\lambda)^b\right).
$$

Thus, we obtain

\begin{prop}\label{gauss} If $\lambda\neq 0$, the number of ${\mathbb F}_q$-rational points in
$X_\lambda$ is given by
$$
N_q(\lambda)=\frac{(q-1)^n-(-1)^n}{q}+\frac{(-1)^n}{q-1}+\frac{(-1)^n}{q(q-1)}\sum_{(n+1)a+b\equiv
0 (q-1) \atop (a,b)\neq(0,0)} {G_q(a)^{n+1}
G_q(b)}\omega(-\lambda)^b.
$$
\end{prop}

We can rewrite this as
$$
N_q(\lambda)=\frac{(q-1)^n-(-1)^n}{q}+\frac{(-1)^n}{q-1}+\frac{(-1)^n}{q(q-1)}\sum_{k=1}^{q-2}
{G_q(k)^{n+1} G_q(-(n+1)k)}\omega(-\lambda)^{-(n+1)k}.
$$

If $\lambda=0$, then equation \ref{inv} gives
$$
S_q(0)=\sum_{x_0,x_1,\ldots,x_n\in\Fq^\star}\zeta_p^{\Tr(x_0x_1+\cdots+x_0x_n+\frac{x_0}{x_1\cdots
x_n})}$$
$$
=\sum_{x_0,x_1,\ldots,x_n\in\Fq^\star}\zeta_p^{\Tr(x_0x_1)}\cdots\zeta_p^{\Tr(x_0x_n)}\zeta_p^{\Tr(x_0/x_1\cdots
x_n)}$$
$$
=\sum_{k_1,\ldots,k_{n+1}=0}^{q-2}\frac{G_q(k_1)\cdots
G_q(k_{n+1})}{(1-q)^{n+1}}\sum_{y_i^{q-1}=1}(y_0y_1)^{k_1}\cdots
(y_0y_n)^{k_n}(\frac{y_0}{y_1\cdots y_n})^{k_{n+1}}
$$
$$
=\sum_{k_1,\ldots,k_{n+1}=0}^{q-2}\frac{G_q(k_1)\cdots
G_q(k_{n+1})}{(1-q)^{n+1}}\sum_{y_i^{q-1}=1}y_0^{k_1+\cdots+k_{n+1}}y_1^{k_1-k_{n+1}}\cdots
 y_n^{k_n-k_{n+1}}
$$
$$
=(-1)^{n+1}\sum_{{k=0} \atop (n+1)k\equiv 0 (q-1)}^{q-2}{G_q(k)^{n+1}}=(-1)^{n+1}\left(1+\sum_{{k=1} \atop (n+1)k\equiv 0 (q-1)}^{q-2}{G_q(k)^{n+1}}\right).
$$

And therefore
$$
N_q(0)=\frac{(q-1)^n-(-1)^n}{q}+\frac{(-1)^{n+1}}{q}\sum_{{k=1} \atop (n+1)k\equiv 0
(q-1)}^{q-2}{G_q(k)^{n+1}}.
$$
Writing $(n+1)=p^am$ with $(p,m)=1$, this is
\begin{equation}\label{number}
N_q(0)=\frac{(q-1)^n-(-1)^n}{q}+\frac{(-1)^{n+1}}{q}\sum_{{k=1} \atop mk\equiv 0
(q-1)}^{q-2}{G_q(k)^{n+1}}.
\end{equation}

Let $S_m=\{\frac{1}{m},\frac{2}{m},\cdots,\frac{m-1}{m}\}$. It is clear that
multiplication by $p$ induces an action on $S_m$, called $p$-action:
$$
r\mapsto \{pr\}
$$
where $\{pr\}$ denotes the fractional part of $pr$. For a given
$r\in S_m$, let $d(r)$ denote the length of the $p$-orbit
containing $r$, that is, the smallest positive integer $d$ such
that $(p^d-1)r\in{\mathbb Z}$. Let $S_{m,d}$ denote the set of
$p$-orbits of length $d$ in $S_m$.

Since $G_{p^d}(k)=G_{p^d}(pk)$, it is clear that if $r_1$ and
$r_2$ are in the same $p$-orbit $\sigma$ in $S_{m,d}$,
$G_{p^d}(r_1(p^d-1))=G_{p^d}(r_2(p^d-1))$. Let us denote this
common value by $G_{p^d}(\sigma (p^d-1))$. Since the set of
$p$-orbits of $S_m$ is the union of $S_{m,d}$ for all $d\geq 1$,
we have

\begin{thm}\label{zetazero} The zeta function of $X_0$ over ${\mathbb F}_p$ is
given by
$$
Z(X_0,T)^{(-1)^n}=\prod_{i=0}^{n-1}(1-p^{i}T)^{{n\choose
i+1}(-1)^{i}}\prod_{d\geq 1}\prod_{\sigma\in
S_{m,d}}\left(1-T^d\frac{G_{p^d}^{n+1}(\sigma(p^d-1))}{p^d}\right).
$$
\end{thm}

{\bf Proof.} By equation \ref{number},
$$
\log Z(X_0,T)=\sum_{k\geq 1} \frac{T^k}{k}\frac{(p^k-1)^n-(-1)^n}{p^k}+\sum_{k\geq
1}\frac{T^k}{k}\frac{(-1)^{n+1}}{p^k}\sum_{{h=1} \atop mh\equiv 0
(p^k-1)}^{p^k-2}{G_{p^k}(h)^{n+1}}.
$$
The second sum is
$$
\sum_{k\geq 1}\frac{T^k}{k}\frac{(-1)^{n+1}}{p^k}\sum_{{r\in S_m} \atop
{r(p^k-1)\in{\mathbb Z}}}{G_{p^k}}(r(p^k-1))^{n+1}$$
$$
=(-1)^{n+1}\sum_{d\geq 1}\sum_{\sigma\in
S_{m,d}}\sum_{r\in\sigma}\left(\sum_{k\geq 1}
\frac{T^{dk}}{dk}\frac{{G_{p^{dk}}}(r(p^{dk}-1))^{n+1}}{p^{dk}}\right)$$
$$
=(-1)^{n+1}\sum_{d\geq 1}\sum_{\sigma\in S_{m,d}}\left(\sum_{k\geq 1}
\frac{T^{dk}}{k}\frac{{G_{p^{dk}}}(\sigma(p^{dk}-1))^{n+1}}{p^{dk}}\right).
$$
By the Hasse-Davenport relation, this sum becomes
$$
=(-1)^{n+1}\sum_{d\geq 1}\sum_{\sigma\in S_{m,d}}\left(\sum_{k\geq 1}
\frac{T^{dk}}{k}\frac{G_{p^d}(\sigma(p^{d}-1))^{k(n+1)}}{p^{dk}}\right),
$$
which gives the stated formula for the zeta function.\hfill$\Box$

\begin{cor} 1) The rank of the sheaf $\mathcal F$ at $0$ is $m-1$.

2) The local $L$-function of the sheaf $\mathcal F$ at $0$ is given by
$$
\prod_{d\geq 1}\prod_{\sigma\in
S_{m,d}}\left(1-T^d\frac{G_{p^d}^{n+1}(\sigma(p^d-1))}{p^d}\right).
$$
\end{cor}

{\bf Proof.} From the given formula for the $L$-function, we see that the degree of
the non-trivial part is given by $|S_m|=m-1$.\hfill$\Box$.


\begin{thebibliography}{}


\bibitem{bbd} Beilinson, A., Bernstein, J. and Deligne, P.
\newblock \emph{Faisceaux pervers,} \newblock in \emph{Analyse et
topologie sur les espaces singuliers (I).} \newblock Ast\'erisque
100. \newblock Societ\'e Math\'ematique de France (1982).

\bibitem{Bu} Buzzard, K.
\newblock \emph{Questions on slopes of modular forms,}
\newblock Ast\'erisque 298 (2005), 1-15.

\bibitem{Can1} Candelas, P., de la Ossa, X., Rodriques-Villegas, F.,  \newblock \emph{Calabi-Yau manifolds over
finite fields, I.}
\newblock http://xxx.lanl.gov/abs/hep-th/0012233.

\bibitem{Can2} Candelas, P., de la Ossa, X., Rodriques-Villegas, F.,  \newblock \emph{Calabi-Yau manifolds over
finite fields, II.}
\newblock Fileds Inst. Commun., \textbf{38}(2003).

\bibitem{CM} Coleman.R and Mazur. B.
\newblock \emph{The Eigencurve,} \newblock in \emph{Galois Representations
in Arithmetic Geometry}(Durham, 1996), \newblock 1-113, London
Math. Soc. Lecture Note Ser., 254,
\newblock Cambridge Univ. Press. (1998).


\bibitem{weilii} Deligne, P. \newblock \emph{La Conjecture de Weil
II.}
\newblock Inst. Hautes \'Etudes Sci. Publ. Math. \textbf{52} (1980)
137--252.

\bibitem{deligne} Deligne, P. \newblock \emph{Application de la
formule des traces aux sommes trigonom\'etriques,} \newblock in
\emph{Cohomologie \'Etale (SGA 4 1/2)} \newblock (Springer-Verlag
1977), 168--232.

\bibitem{sga7ii} Deligne, P. and Katz, N. \newblock \emph{Groupes de Monodromie en G\'eom\'etrie
Alg\'ebrique (SGA 7 II).} \newblock Lecture Notes in Mathematics
340, Springer-Verlag (1973).


\bibitem{Dw0} Dwork, B. \newblock \emph{$p$-adic Cycles.}
\newblock Publ. Math. IHES,  \textbf{37} (1969), 27-115.

\bibitem{Dw} Dwork, B. \newblock \emph{Normalized period matrices II.}
\newblock Ann. Math., \textbf{98} (1973), 1-57.

\bibitem{fw0} Fu, L. and Wan, D. \newblock \emph{Moment zeta functions, partial L-functions and partial
exponential sums.}
\newblock Math. Ann., \textbf{328} (2004), 193-228.

\bibitem{fw} Fu, L. and Wan, D. \newblock \emph{L-functions for
symmetric products of Kloosterman sums.} \newblock J. Reine Angew.
Math. \textbf{589} (2005), 79--103.

\bibitem{fw2} Fu, L. and Wan, D. \newblock \emph{Trivial factors for L-functions of
symmetric products of Kloosterman sheaves.} arXiv:math.AG/0610228
v1 6 Oct 2006.

\bibitem{fw3} Fu, L. and Wan, D. \newblock \emph{Mirror congruence for rational points
on Calabi-Yau varieties.} \newblock Asian J.
Math. \textbf{10} (2006), no. 1, 1-10.

\bibitem{fh} Fulton, W and Harris, J. \newblock \emph{Representation
Theory.} \newblock Graduate Texts in Mathematics, 129. \newblock
Springer-Verlag (1991).

\bibitem{Ta} Harris, M., Shepherd-Barron, N. and Taylor, R.
\newblock \emph{A family of Calabi-Yau varieties and potential automorphy.}
\newblock Preprint, 2006.

\bibitem{ka0} Katz, N. \newblock \emph{Algebraic solutions of differential equations ($p$-curvature
and the Hodge filtration.} \newblock Invent. Math.,
\textbf{18} (1972), 1-118.

\bibitem{ka1} Katz, N. \newblock \emph{Sommes Exponentielles.} \newblock Asterisque
\textbf{79}(1980), 209pp.

\bibitem{gkm} Katz, N. \newblock \emph{Gauss Sums, Kloosterman
Sums and Monodromy Groups.} \newblock Annals of Mathematics Studies
116 (Princeton University Press 1987).

\bibitem{esde} Katz, N. \newblock \emph{Exponential Sums and Differential Equations.} \newblock Annals of Mathematics Studies
124 (Princeton University Press 1990).

\bibitem{fs} Katz, N. \newblock \emph{Frobenius-Schur Indicator and
the Ubiquity of Brock-Granville Quadratic Excess.} \newblock Finite Fields Appl.
\textbf{7} (2001), 45--69.

\bibitem{kl} Katz, N. and Laumon, G. \newblock
\emph{Transformation de Fourier et Majoration de Sommes
Exponentielles.}
\newblock   Inst. Hautes \'Etudes Sci. Publ. Math. \textbf{62} (1985)
361--418.

\bibitem{k2} Katz, N., \newblock Another look at the Dwork family, \newblock preprint, 2007.


\bibitem{laumon} Laumon, G. \newblock \emph{Transformation de
Fourier, constantes d'\'equations fonctionelles et conjecture de
Weil.} \newblock Pub. Math. IHES \textbf{65} (1987), 131--210.

\bibitem{ko} Koblitz, N. \newblock \emph{The number of points on certain family of hypersurfaces over finite fields.}
\newblock Compositio Math., \textbf{48} (1983), no.1, 3-23.


\bibitem{W0} Wan, D. \newblock \emph{Partial zeta functions of algebraic varieties over finite fields.}
\newblock Finite Fields \& Applications, \textbf{7} (2001), 238-251.

\bibitem{W1} Wan, D. \newblock \emph{Dimension variation of classical and $p$-adic modular forms.}
\newblock Invent. Math., \textbf{133} (1998), 449-463.

\bibitem{W20} Wan, D. \newblock \emph{Dwork's conjecture on unit root zeta functions.}
\newblock Ann. Math., \textbf{150} (1999), 867-927.

\bibitem{W21} Wan, D. \newblock \emph{Higher rank case of Dwork's conjecture.}
\newblock J. Amer. Math. Soc., \textbf{13} (2000), 807-852.

\bibitem{W2} Wan, D. \newblock \emph{Rank one case of Dwork's conjecture.}
\newblock J. Amer. Math. Soc., \textbf{13} (2000), 853-908.

\bibitem{W3} Wan, D. \newblock \emph{Geometric moment zeta functions.}
\newblock In Geometric Aspects of Dwork Theory, Walter de Gruyter, 2004,
VolII, 1113-1129.

\bibitem{W31} Wan, D. \newblock \emph{Variation of $p$-adic Newton polygons
for L-functions of exponential sums.}
\newblock Asian J. Math., \textbf{8} (2004), 427-474.

\bibitem{W4} Wan, D. \newblock \emph{Arithmetic mirror symmetry.}
\newblock Pure Appl. Math. Q., Vol. 1(2005) no. 2, 369-378.

\bibitem{W5} Wan, D. \newblock \emph{Mirror symmetry for zeta functions.}
\newblock In Mirror Symmetry V, AMS/IP Studies in Advanced Mathematics, Vol. 38,
2007.


\end{thebibliography}
\end{document}